\documentclass[9pt]{article}

\usepackage{amsmath, amsthm}
\usepackage{amsmath, amsfonts}
\usepackage{amsmath, amssymb}
\usepackage{amsmath}
\usepackage[all]{xy}

\newtheorem{theorem}{Theorem}[section]
\newtheorem{definition}[theorem]{Definition}
\newtheorem{definition-lemma}[theorem]{Definition/Lemma}
\newtheorem{definition-explanation}[theorem]{Definition/Explanation}
\newtheorem{explanation-definition}[theorem]{Explanation/Definition}
\newtheorem{definition-fact}[theorem]{Definition/Fact}
\newtheorem{definition-notation}[theorem]{Definition/Notation}
\newtheorem{proto-definition}[theorem]{Proto-Definition}
\newtheorem{quasi-definition}[theorem]{Quasi-Definition}
\newtheorem{lemma}[theorem]{Lemma}
\newtheorem{lemma-definition}[theorem]{Lemma/Definition}
\newtheorem{proposition}[theorem]{Proposition}

\newtheorem{remark}[theorem]{\it Remark}
\newtheorem{remark-notation}[theorem]{\it Remark/Notation}
\newtheorem{conjecture}[theorem]{Conjecture}

\newtheorem{example}[theorem]{Example}
\newtheorem{example-definition}[theorem]{Example/Definition}

\newtheorem{definition-prototype}[theorem]{Definition-Prototype}

\numberwithin{equation}{section}

\input{epsfig.sty}

 \newcommand{\Azscriptsize}{\mbox{\scriptsize\it A$\!$z}}

\newcommand{\End}{\mbox{\it End}\,}
 \newcommand{\smallEnd}{\mbox{\small\it End}\,}
\newcommand{\Endsheaf}{\mbox{\it ${\cal E}\!$nd}\,}

\newcommand{\GL}{\mbox{\it GL}}

\newcommand{\Hom}{\mbox{\it Hom}\,}

\newcommand{\Image}{\mbox{\it Im}\,}

\newcommand{\Ker}{\mbox{\it Ker}\,}

\newcommand{\Mod}{\mbox{\it Mod}\,}

\newcommand{\Quot}{\mbox{\it Quot}}

\newcommand{\Space}{\mbox{\it Space}\,}

\newcommand{\Spec}{\mbox{\it Spec}\,}
 \newcommand{\boldSpec}{\mbox{\it\bf Spec}\,}
 
 \newcommand{\smallSpec}{\mbox{\small\it Spec}\,}

\newcommand{\Supp}{\mbox{\it Supp}\,}
 
\newcommand{\Sym}{\mbox{\it Sym}}

\newcommand{\bii}{\mbox{\it b.i.i.}}

\newcommand{\et}{\mbox{\scriptsize\it \'{e}t}\,}

\newcommand{\pt}{\mbox{\it pt}}

\begin{document}

\noindent

\centerline{\large\bf
 Azumaya noncommutative geometry and D-branes}
\vspace{1ex}
\centerline{\large\bf
 - an origin of the master nature of D-branes}
\vspace{2em}
\centerline{\large Chien-Hao Liu}
\vspace{2em}

\baselineskip 10pt
\begin{quotation} {\small
{\bf Abstract.}
 In this lecture I review how a matrix/Azumaya-type noncommutative
  geometry arises for D-branes in string theory and
  how such a geometry serves as an origin of the master nature of D-branes;
  and then highlight an abundance conjecture on D0-brane resolutions
  of singularities that is extracted and purified
  from a work of Douglas and Moore in 1996.
 A conjectural relation of our setting with `D-geometry' in the sense of
  Douglas is also given.
 The lecture is based on a series of works on D-branes with Shing-Tung Yau,
  and in part with Si Li and Ruifang Song.
}\end{quotation}
{\small
 Parts delivered in the workshop
  {\sl Noncommutative algebraic geometry and D-branes},
  December 12 -- 16, 2011,
  organized by Charlie Beil, Michael Douglas, and Peng Gao,
  at Simons Center for Geometry and Physics,
  Stony Brook University, Stony Brook, NY.}
\vspace{1em}

\noindent{\small
{\sc Dedication.}
 This lecture is dedicated to
     {\it Shiraz Minwalla},
 $\,${\it Mihnea Popa},
 $\,${\it Ling-Miao Chou},
 who together made this project possible;
 and to my mentors (time-ordered): {\it Hai-Chau Chang},
 {\it William Thurston}, {\it Orlando Alvarez},
 {\it Philip Candelas}, {\it Shing-Tung Yau},
 who together shaped my unexpected stringy/brany path.}

\vspace{8em}

\begin{flushleft}
{\small\bf Outline.}
\end{flushleft}
{\small
\begin{itemize}
 \item[1.]
  D-brane as a morphism from Azumaya noncommutative spaces
  with a fundamental module.
  \vspace{-.6ex}
  \begin{itemize}
   \item[$\cdot$]
    The emergence of a matrix-/Azumaya-type noncommutativity.

   \item[$\cdot$]
    A naive/direct space-time interpretation of this noncommutativity.

   \item[$\cdot$]
    A second look:
    What is a D-brane (mathematically)? - From Polchinski to Grothendieck.

   \item[$\cdot$]
    What is a noncommutative (algebraic) geometry? -
    Looking for a D-brane-sensible/motivated\\ settlement
    in an inperfect noncommutative world.

   \item[$\cdot$]
    Reflection and a conjecture on D-geometry in the sense of Douglas:\\
    Douglas meeting Polchinski-Grothendieck.
  \end{itemize}

 \item[2.]
  Azumaya geometry as the origin of the master nature of D-branes.
  \vspace{-.6ex}
  \begin{itemize}
   \item[$\cdot$]
    Azumaya noncommutative geometry
    as the origin of the master nature of D-branes.

   \item[$\cdot$]
    Azumaya noncommutative algebraic geometry
    as the master geometry for commutative algebraic geometry.
  \end{itemize}

 \item[3.]
  D-brane resolution of singularities - an abundance conjecture.
  \vspace{-.6ex}
  \begin{itemize}
   \item[$\cdot$]
    Beginning with Douglas and Moore:
    D-brane resolution of singularities.

   \item[$\cdot$]
    The richness and complexity of Azumaya noncommutative space.

   \item[$\cdot$]
    An abundance conjecture.
  \end{itemize}

 \item[] \hspace{-1.6em}
  Epilogue.\\[.2ex]
  $\mbox{\hspace{-1.6em}}$
  Notes and acknowledgements added after the workshop.\\[.2ex]
  $\mbox{\hspace{-1.6em}}$
  References.
\end{itemize}
} 

\newpage
\baselineskip 12pt

\vspace{2em}

\section{D-brane as a morphism from Azumaya noncommutative spaces
         with a fundamental module.}

My lecture today is based on three guiding questions:
 \marginpar{\raggedright\tiny $\bullet$ Prepared on {\bf blackboard}.}
 \begin{itemize}
  \item[]
   {\bf Q.1} {\it What is a D-brane?}\\
   {\bf Q.2} {\it What is a noncommutative geometry?}\\
   {\bf Q.3} {\it How are the two related?}
 \end{itemize}
To reflect the background of this lecture, I assume:
 \marginpar{\raggedright\tiny $\bullet$ Prepared on {\bf blackboard}.}
 \begin{itemize}
  \item[]
   {\bf When:} October, 1995; or, indeed, 1989.\\
   {\bf Where:}
    In the {\it geometric phase of} Wilson's theory-space
     ${\cal S}_{Wilson}^{d=2, CFT w/ boundary}$
     for $d=2$ conformal field theory with boundary;$\,/\!/\hspace{1ex}$
    and with assumption that {\it open string tension is large enough}
     (so that {\it D-brane is soft} with respect to open strings).
 \end{itemize}

\bigskip

\begin{flushleft}
{\bf The emergence of a matrix-/Azumaya-type noncommutativity.}
\end{flushleft}
$\bullet$
 Let me begin with Polchinski's TASI lecture on D-branes in 1996 ...
 \begin{itemize}
  \item[$\cdot$]
   ... and first recall the very definition of a D-brane
   from string theory:

   \begin{definition} {\bf [D-brane].} {\rm
    A {\it Dirichlet-brane} (in brief {\it D-brane})
     is a submanifold/cycle/locus in an open-string target space-time
     in which the boundary/end-points of an open string can lie.}
   \end{definition}

  \item[$\cdot$]
   {\it Figure~1-1:}
   {\it Oriented open strings with end-points on D-branes}.
    \marginpar{\raggedright\tiny $\bullet$ Color chalks.}
    \begin{itemize}
     \item[-]
      $f:X\rightarrow Y$,
      where
       $X$ is endowed with local coordinates $\xi:=(\xi^a)_a$,
       $Y$ local coordinates $(y^a; y^{\mu})_{a,\mu}$, and
       $f$ is given by $y^a=\xi^a$ and $y^{\mu}=f^{\mu}(\xi)$.
    \end{itemize}

  \item[$\cdot$]
   This definition,
    though mathematically far from obvious at all as what it'll lead to,
    is very fundamental from physics point of view.$\,/\!/\hspace{1ex}$
   It says that all the fields on D-branes and the dynamical law
    that governs them are created by open strings.

  \item[$\cdot$]
   Open strings vibrate and its end-points create
    (both massless and massive) fields
     on the D-brane world-volume.$\,/\!/\hspace{1ex}$
   Massless fields are created by an open string
    with both ends on the same branes.$\,/\!/\hspace{1ex}$
   There are two complementary sets of these:
    One corresponds to vibrations of ends of the open string
     in the tangential directions along the D-brane.
    This creates an $u(1)$ gauge field on the branes.
    The other set corresponds to vibrations of ends of the open string
        in the normal directions to the D-brane.
    This creates a scalar field that describes fluctuations of the D-brane
     in space-time.
 \end{itemize}

\bigskip

\noindent $\bullet$
When $r$-many D-branes coincide in space-time,
 something mysterious happens:
 \begin{itemize}
  \item[$\cdot$]
   One key feature of an open or closed string,
     compared to the usual mechanical string in our daily life,
    is that its tension is a constant in the theory;$\,/\!/\hspace{1ex}$
   and hence the mass of states or fields on D-branes created
    by open-strings
    are proportional to the length of the string.$\,/\!/\hspace{1ex}$
   Once $r$-many D-branes are brought to coincide in space-time,
    there are states/fields that were originally massive
    but now becomes massless.$\,/\!/\hspace{1ex}$
   (Continuing {\it Figure~1-1}.)

  \item[$\cdot$]
   In particular, the gauge fields $A_{a}$ on the stacked D-brane
    is now enhanced to $u(r)$-valued$\,/\!/\hspace{1ex}$
   and the scalar field $y^{\mu}$ on the D-brane world-volume
   that describes the deformation of the brane is also $u(r)$-valued.

  \item[$\cdot$]
   For this, Polchinski made the following comment in his by-now-standard
   textbook for string theory:
   \begin{itemize}
    \item[$\cdot$] ([Po2: vol.~I, Sec.~8.7, p. 272].)
    \marginpar{\raggedright\tiny $\bullet$
      Prepared on {\bf blackboard}.}
    (With mild notation change.)\\
    ``For $r$-separated D-branes, the action is $r$ copies of the action
       for a single D-brane. We have seen, however, that {\it
       when the D-branes are coincident, there are $r^2$ rather than $r$
       massless} vectors and {\it scalars on the brane}, and we would like
       to write down the effective action governing these.
      The fields $y^{\mu}(\xi)$ and $A_a(\xi)$
       will now be $r\times r$ matrices.
      For the gauge field, the meaning is obvious -- it becomes a non-Abelian
       $U(r)$ gauge field.
      {\it For the collective coordinates $y^{\mu}$, however,
        the meaning is mysterious:
       the collective coordinates for the embedding of $r$ D-branes
        in spacetime are now enlarged to $r\times r$ matrices.}
      This `{\it noncommutative geometry}' has proven to play a key role
       in the dynamics of D-branes, and there are conjectures that
       it is an important {\it hint about the nature of spacetime}."
   \end{itemize}
 \end{itemize}

\bigskip

\begin{flushleft}
{\bf A naive/direct space-time interpretation of this noncommutativity.}
\end{flushleft}
$\bullet$
As $y^{\mu}$ are meant to be the coordinates for the open-string
 target-space-time $Y$,
it is very natural for one to perceive that
 somehow there is something noncommutative about this space-time
 that is originally hidden from us
 before we let the D-branes collide.$\,/\!/\hspace{1ex}$
And once we let the D-branes collide,
 this hidden feature of space-time reveals itself suddenly
 through a new geometry
 whose coordinates are matrix/Azumaya-algebra-valued.$\,/\!/\hspace{1ex}$
It seems to me that this is what Polchinski reflects in the above comment
 and it turns out to be what the majority of stringy community
 think about as well.

\bigskip

\begin{flushleft}
{\bf A second look:
     What is a D-brane (mathematically)? - From Polchinski to Grothendieck.}
\end{flushleft}
$\bullet$
Re-think about the phenomenon
 {\it locally} and {\it from Grothendieck's construction}
  of modern algebraic geometry via the language schemes:
\begin{itemize}
 \item[$\cdot$]
  Let $R(X)$ be the ring of local functions
   (e.g.\ $C^{\infty}(X)$ in real smooth category)
   of $X$ and $R(Y)$ be the ring of local functions on $Y$
   (e.g.\ $C^{\infty}(Y)$).$\,/\!/\hspace{1ex}$
  Then $\xi^a\in R(X)\,$; $\,y^a,\, y^{\mu}\in R(Y)\,$; and
   $f$ above is equivalently but contravariantly
   given by a ring-homomorphism
   $f^{\sharp}: R(Y)\rightarrow R(X)$
   specified by
   $$
    y^a\longmapsto \xi^a
     \hspace{2em}\mbox{and}\hspace{2em}
    y^{\mu}\longmapsto f^{\mu}(\xi)\,,
   $$
   i.e. $f:X\rightarrow Y$ is determined
    how it {\it pulls back local functions from $Y$ to $X$}.

 \item[$\cdot$]
  When $r$-many D-branes coincide, formally $y^{\mu}$ becomes matrix-valued.
  But $y^{\mu}$ takes values in the function ring of $X$
   under $f^{\sharp}$.$\,/\!/\hspace{1ex}$
  This suggest that the original $R(X)$ is now enhanced to
   $M_r(R(X))$
   (or more precisely
       $M_r(R(X)\otimes_{\Bbb R}{\Bbb C})
        =M_r({\Bbb C})\otimes_{\Bbb R}R(X)$).$\,/\!/\hspace{1ex}$
  In other words,
   {\it the D-brane world-volume becomes matrix/Azumaya noncommutatized!}
\end{itemize}

\begin{remark} {$[\,$pure open-string effect$\,]$.} {\rm
 It is conceptually worth emphasizing that, from the above reasoning,
  one deduces also that
  {\it this fundamental noncommutativity on D-brane world-volume
   is a purely open-string induced effect.}$\,/\!/\hspace{1ex}$
 No $B$-field, supersymmetry, or any kind of quantization is involved.
}\end{remark}

\begin{remark}
{$[\,$Lie algebra vs.\ Azumaya/matrix-ring algebra$\,]$.} {\rm
 Acute string theorists may recall that in the original string-theory
  setting and in the world-volume field-theory language,
  this field $y^{\mu}$ is indeed an $u(r)$-adjoint scalar.
 So, {\it why didn't we take directly the Lie-algebra-enhancement
  $u(r)\otimes R(X)$ to the function ring $R(X)$ of
  the D-brane world-volume $X$?}$\,/\!/\hspace{1ex}$
 The answer comes from two sources:
  \begin{itemize}
   \item[(1)] {\it For geometry reason}$\,:$
    Local function ring of a geometry has better to be associate
     and with an identity element $1$.$\,/\!/\hspace{1ex}$
    Without the latter, one doesn't even know how to start
     for a notion of localization of the ring,
     a concept that is needed for a local-to-global gluing construction.

   \item[(2)] {\it For field-theory reason}$\,:$
    The {\it kinetic term} is the action on D-brane world-volume
    involves matrix multiplication;
    it is not expressible in terms of Lie brackets alone.
  \end{itemize}
}\end{remark}

\begin{proto-definition} {\bf [D-brane: Polchinski-Grothendieck].} {\rm
 A {\it D-brane} is
  {\it an Azumaya noncommutative space with a fundamental module}
  $$
   (X^{A\!z},{\cal E}):=(X,{\cal O}_X^{A\!z},{\cal E})\,,
  $$
   where ${\cal O}_X^{A\!z}=\Endsheaf_{{\cal O}_X}({\cal E})$.
 A {\it D-brane on $Y$} is a morphism
  $$
   \varphi\; :\; (X^{A\!z},{\cal E})\; \longrightarrow\; Y
  $$
  {\it defined by}
  $$
   \varphi^{\sharp}\; :\;
    {\cal O}_Y\; \longrightarrow\; {\cal O}_X^{A\!z}
  $$
  {\it as an equivalence class of gluing systems of ring homomorphisms
  of local function rings from $Y$ to $X$}.
}\end{proto-definition}

\noindent $\bullet$
Two reasons I call this a proto-definition for D-branes:
 \begin{itemize}
  \item[(1)]
   I {\it focus only on fields on D-branes that are relevant to
    the occurrence of the matrix/Azumaya type noncommutativity in question}.

  \item[(2)]
   I conceal {\it subtle local-to-global issues from the constructibility
    and nonconstructibility in noncommutative geometry},
   which I need to explain and will come back ...
 \end{itemize}
... but, to help casting away the possible doubt from string theorists
 as whether this proto-definition makes sense,
 let me give first a very simple, concrete, and yet deep enough example
 which we are now ready.

\begin{example}
{\bf [D0-brane on the complex line ${\Bbb A}^1_{\Bbb C}$
      via Polchinki-Grothendieck].} {\rm
 An {\it Azumaya point$/{\Bbb C}$ with a fundamental module of rank $r$}
  is given by
  $$
   (pt, \End_{\Bbb C}(E), E)\,,
  $$
  where $E$ is isomorphic to ${\Bbb C}^r$.
 This is our {\it D0-brane}.$\,/\!/\hspace{1ex}$
 To be explicit, let's fix an isomorphism $E\simeq {\Bbb C}^r$,
  which fixes also the ${\Bbb C}$-algebra isomorphism
  $\End_{\Bbb C}(E)\simeq$ the ${\Bbb C}$-algebra $M_r({\Bbb C})$
   of $r\times r$ matrices.
 One should think of this
  as a noncommutative point
  $$
   \Space(M_r({\Bbb C}))\,,
  $$
  whose function ring is given by $M_r({\Bbb C})$,
  with a built-in module ${\Bbb C}^r$
                         of the function ring.$\,/\!/\hspace{1ex}$
 We take the complex line ${\Bbb A}^1_{\Bbb C}$
  as an affine variety over ${\Bbb C}$,
  whose local rings is given the polynomial ring ${\Bbb C}[y]$
  over ${\Bbb C}$ in one variables $y$.
 One could think of this $y$ as a {\it coordinate function}
  on ${\Bbb A}^1_{\Bbb C}.$\,/\!/\hspace{1ex}
 In algebro-geometric notation (and with a few subtleties concealed),
  $$
   {\Bbb A}^1_{\Bbb C}\; =\; \Spec({\Bbb C}[y])\,.
  $$

 Following the setting above,
  a {\it D0-brane on ${\Bbb A}^1_{\Bbb C}$}
   is then a {\it morphism}
   $$
    \varphi\; :\; (\Space(M_r(\Bbb C)),{\Bbb C}^r)\;
                  \longrightarrow\; {\Bbb A}^1_{\Bbb C}
   $$
  {\it defined by a ${\Bbb C}$-algebra homomorphism}
   $$
    \varphi^{\sharp}\;:\;
     {\Bbb C}[y]\; \longrightarrow\; M_r({\Bbb C})\,.
   $$
 This, in turn, is determined by an (arbitrary) specification
   $$
    y\;\longmapsto\; m_{\varphi}\,\in\,M_r({\Bbb C})\,.
   $$

 Now comes the most essential question:
  \begin{itemize}
   \item[{\bf Q.}]
    {\it Does this match with
         how D-branes behave in string-theorists' mind?}
  \end{itemize}
 Let's now examine this by looking at two things:
  \begin{itemize}
   \item[(1)]
    the image 0-brane with Chan-Paton sheaf on ${\Bbb A}^1_{\Bbb C}$;

   \item[(2)]
    how do they vary when we vary $\varphi$.
  \end{itemize}
 Here, we adopt the standard set-up of Grothendieck's theory
  of (commutative) schemes:
  \begin{itemize}
   \item[(1)]
    The image $0$-brane $\Image\varphi$ on ${\Bbb A}^1_{\Bbb C}\,$:
    \begin{itemize}
     \item[-]
      This is the subscheme of ${\Bbb A}^1_{\Bbb C}$ defined by the ideal
       $I_{\varphi} := \Ker\varphi^{\sharp}
                     =(\varphi^{\sharp})^{-1}(0) \subset {\Bbb C}[y]$.

     \item[-]
      Let $I_{\varphi}=((y-c_1)^{n_1}\,\cdots\,(y-c_k)^{n_k})$.
      Then
       $(y-c_1)^{n_1}\,\cdots\,(y-c_k)^{n_k}$
       is the minimal polynomial for $m_{\varphi}$.
      In particular,
       $n_1+\,\cdots\,n_k\le r$ and, ignoring multiplicity,
       $\{c_1,\,\cdots\,,c_k\}$
        is exactly the set of eigen-values of $m_{\varphi}$.

     \item[-]
      In plain words, this says that
       $\Image\varphi$ is a collection of fuzzy/thick points supported
       at points $c_1,\,\cdots,\,c_k$ in the complex line ${\Bbb C}$
       with multiplicity of fuzziness $n_1,\,\cdots,\,n_k$ respectively.
    \end{itemize}

   \item[$\cdot$]
    The Chan-Paton sheaf $\varphi_{\ast}({\Bbb C}^r)\,$:
    \begin{itemize}
     \item[-]
      Through the ${\Bbb C}$-algebra homomorphism
      $\varphi^{\sharp}:{\Bbb C}[y]\rightarrow M_r({\Bbb C})$,
      the $M_r({\Bbb C})$-module ${\Bbb C}^r$ becomes a ${\Bbb C}[y]$-module
       with $I_{\varphi}\cdot {\Bbb C}^r=0$.$\,/\!/\hspace{1ex}$
      Thus, $\varphi_{\ast}({\Bbb C}^r)$ is simply
       ${\Bbb C}^r$ as a ${\Bbb C}[y]/I_{\varphi}$-module.

     \item[-]
      Geometrically, this says that
       $\varphi_{\ast}({\Bbb C}^r)$ is a $0$-dimensional coherent sheaf
       on ${\Bbb A}^1_{\Bbb C}$, supported on the 0-dimensional subscheme
       $\Image\varphi$ of ${\Bbb A}^1_{\Bbb C}$.
    \end{itemize}

   \item[(2)]
    Deformations of $\varphi$ are {\it defined by} deformations
     of the ${\Bbb C}$-algebra homomorphism
     $\varphi^{\sharp}$.$\,/\!/\hspace{1ex}$ \\
    The corresponding $\Image\varphi$ and $\varphi_{\ast}({\Bbb C}^r)$
     on ${\Bbb A}^1_{\Bbb C}$ vary accordingly.
  \end{itemize}
  These are illustrated in {\sc Figure}~1-2.
  \marginpar{\raggedright\tiny $\bullet$ Prepared on {\bf blackboard}.}
 From this very explicit example/illustration, we see that:
  \begin{itemize}
   \item[$\cdot$]
    The notion of {\it Higgsing and un-Higgsing} of D-branes and
     of {\it recombinations} of D-branes
     are nothing but outcomes of {\it deformations of morphisms}
     from an Azumaya space with a fundamental module,
     as is defined in Proto-Definition~1.4.
  \end{itemize}
 In other words, our setting does indeed capture some key features
  of D-branes in string theory!

\noindent\hspace{15cm}$\square$
}\end{example}

\begin{remark}
{$[\,$D-brane world-volume vs.\ open-string target-space-time$\,]$.}
{\rm
 Now we have two aspects of this matrix/Azumaya-type noncommutativity:
  one as part of a hidden structure of open-string target-space-time
  revealed through stacked D-branes,
 and the other as a fundamental structure on the D-brane world-volume
  when D-branes become coincident.$\,/\!/\hspace{1ex}$
 There are two fundamental reasons we favor the latter,
  rather than the former:
  \begin{itemize}
   \item[(1)]
    {\it From the physical aspect/a comparison with quantum mechanics}$\,$:
    In quantum mechanics,
     when a particle moving in a space-time
       with spatial coordinates collectively denoted by $x$,
     $x$ becomes operator-valued.$\,/\!/\hspace{1ex}$
    There we don't take the attitude that
     just because $x$ becomes operator-valued,
      the nature of the space-time is changed.$\,/\!/\hspace{1ex}$
    Rather, we say that {\it the particle is quantized but the space-time
     remains classical.}$\,/\!/\hspace{1ex}$
    In other words, it is the nature of the particle that is changed,
     {\it not} the space-time.$\,/\!/\hspace{1ex}$
    Replacing the word `{\it quantized}' by
    `{\it matrix/Azumaya noncommutatized}',
    one concludes that this matrix/Azumaya-noncommutativity happens
     on D-branes, {\it not (immediately on) the space-time}.

   \item[(2)]
    {\it From the mathematical/Grothendieck aspect}$\,$:
    The function ring $R$ is more fundamental
     than the topological space $\Space(R)$, if definable.
    A morphism
     $$
      \varphi\; :\;  \Space(R)\; \longrightarrow\; \Space(S)
     $$
     is specified {\it contravariantly} by a ring-homomorphism
     $$
      \varphi^{\sharp}\; :\; S\; \longrightarrow\; R\,.
     $$
    If the function ring $R$ of the domain space $\Space(R)$
     is commutative, then $\varphi^{\sharp}$ factors through
     a  ring-homomorphism $\bar{\varphi}^{\sharp}:S/[S,S]\rightarrow R$,
     $$
      \xymatrix{
        R  &&  S\ar[ll]_-{\varphi^{\sharp}} \ar@{->>}[d]^-{\pi_{S/[S,S]}}
                \ar@{}[dl]|{\mbox{\raisebox{2.4ex}{$\circ$}}}   \\
           &&  S/[S,S]  \ar[llu]^-{\bar{\varphi}^{\sharp}}   &.
      }
     $$
   Here,
    $[S,S]$, the {\it commutator} of $S$, is the bi-ideal of $S$
     generated by elements of the form $s_1s_2-s_2s_1$
     for some $s_1,\,s_2\in S$;  and
    $S/[S,S]$ is the commutatization of $S$.
   It follows that
     $$
      \xymatrix{
       \Space(R)\ar[rr]^-{\varphi} \ar[rrd]_-{\bar{\varphi}}
        &&  \Space(S)\ar@{}[dl]|{\mbox{\raisebox{2.4ex}{$\circ$}}}\\
        &&  \Space(S/[S,S])\rule{0ex}{2.4ex} \ar@{^{(}->}[u]_-{\iota}   &.
      }
     $$
   In other words,
    \begin{itemize}
     \item[$\cdot$] {\it if
      the function ring on the D-brane world-volume is only commutative,
      then it won't be able to detect the noncommutativity, if any,
      of the open-string target-space!}
    \end{itemize}
  \end{itemize}
  Cf.~{\sc Figure}~1-3.
}\end{remark}

\begin{example}
{$[\;$implicit examples in string theory literature$\,]$.} {\rm
 Once accepting the above aspect
   from Grothendieck's viewpoint of geometry,
  one immediately recognizes that there are many local examples
  hidden implicitly in the string theory literature.
 For instance, the {\it commuting variety/scheme}
  $$
   \{(m_1,\,\cdots\,,\,m_l)\:
          :\; m_i\in M_r({\Bbb C})\,,\;[m_i,m_j]=0\,,\;
                                       1\,\le\,i,j\,\le l\,\}
  $$
  that appears in the description of the D-brane ground states
  in the Coulomb branch/phase of
  the supersymmetric gauge theory coupled with matter
  on the D-brane world-volume
  is exactly the moduli space of morphisms
  from the {\it fixed} Azumaya point-with-a-fundamental module
  $(\Spec{\Bbb C}, M_r({\Bbb C}), {\Bbb C}^r)$ to the affine space
  ${\Bbb A}^l_{\Bbb C}:=\Spec({\Bbb C}[y_1,\,\cdots\,,y_l])$.
 This moduli space in general is quite complicated,
  having many nonreduced irreducible components as a scheme.
 It is indeed canonically isomorphic to the Quot-scheme
  $\Quot({\cal O}_{{\Bbb A}^l_{\Bbb C}}^{\oplus r},r)$ of $0$-dimensional
  coherent ${\cal O}_{{\Bbb A}^l_{\Bbb C}}$-module of length $r$
  on ${\Bbb A}^l_{\Bbb C}$.
 After modding out the global symmetry $\GL_r({\Bbb C})$,
  which corresponds to the change of basis of ${\Bbb C}^r$,
  one obtains the stack
  $$
   {\frak M}^{0^{{A\!z}^f}}\!\!({\Bbb A}^l)\; \simeq\;
   [\Quot({\cal O}_{{\Bbb A}^l_{\Bbb C}}^{\oplus r},r)/\GL_r({\Bbb C})]
  $$
  of D0-branes of length $r$ on ${\Bbb A}^l$.

 For another instance, whenever
  one sees a ring-homomorphism or an {\it algebra representation}
  $$
   \rho\; :\; A\; \longrightarrow\; M_r(B)\,,
  $$
   where
    $A$ is a (possibly noncommutative) associative, unital ring
    -- for example, a quiver algebra --  and
    $B$ is a (usually-commutative-but-not-required-so) ring,
 one is indeed looking at a morphism
  from an Azumaya space with a fundamental module
  $$
   \varphi_{\rho}\; :\;
     (\Space(B), M_r(B), B^{\oplus r})\; \longrightarrow\; \Space(A)
  $$
  {\it defined by} $\rho$, i.e.\ a {\it D-brane on} $\Space(A)\,$!
}\end{example}

\bigskip

\begin{flushleft}
{\bf What is a noncommutative (algebraic) geometry? -
     Looking for a D-brane-sensible/motivated settlement
     in an inperfect noncommutative world.}
\end{flushleft}
$\bullet$
Morphisms between ringed spaces: first attempt.
 \begin{itemize}
  \item[$\cdot$]
   Taking Grothedineck's path:
   (local/affine picture; all rings assumed associative and unital)
   $$
    \mbox{noncommutative ring $R$}\;     \Longrightarrow\;
    \mbox{topological space $\Spec R$}\; \Longrightarrow\;
    \mbox{ringed space $(\Spec R, R)$}\,.
   $$

  \item[$\cdot$]
   A morphism from $(X,{\cal O}_X)\rightarrow (Y,{\cal O}_Y)$
    is given by a pair $(f, f^{\sharp})$,
    where
     $f: X\rightarrow Y$
       is a continuous map between topological spaces  and
     $f^{\sharp}:{\cal O}_Y\rightarrow \varphi_{\ast}{\cal O}_X$
      is a map of sheaves of rings on $Y$.

  \item[$\cdot$]
   Leaving aside the issue of localizations,
   the starting point $R\Rightarrow \Spec R$
    already imposes challenges;
   there are subtle issues on the notion/construction of $\Spec R$
    in the case of general noncommutative rings.
   This remains an ongoing issue for the current and the future
    noncommutative algebraic geometers.
 \end{itemize}

\bigskip

\noindent $\bullet$
Another path via the category of quasi-coherent sheaves.
 \begin{itemize}
  \item[$\cdot$]
   A fundamental work [Ro] of Alexander Rosenberg (1998):
    {\it The spectrum of abelian categories and reconstruction of schemes}.

  \item[$\cdot$]
   Instead of constructing
    noncommutative algebraic geometry from noncommutative rings $R$,
   construct noncommutative geometry from the category $\Mod_R$
    of $R$-modules!

  \item[$\cdot$]
   An unfortunate fact:
   {\it Non-isomorphic noncommutative rings may have
        equivalent categories of modules};
        cf.\ {\it Morita equivalence}.
   That is,
   \begin{itemize}
    \item[$\cdot$]
     {\it in general,
          $\Mod_R$ does not contain all the information of $R$
          when $R$ is noncommutative}.
   \end{itemize}
   Indeed,
   the two ${\Bbb C}$-algebras, $M_r({\Bbb C})$ and ${\Bbb C}$,
    are Morita equivalent.
   More generally:
    \begin{itemize}
     \item[$\cdot$]
      Let $(X,{\cal O}_X)$ be a (commutative) scheme
       and ${\cal E}$ be a locally free sheaf on $X$.
      Then the two sheaves of algebras,
       $\Endsheaf_{{\cal O}_X}({\cal E})$ and ${\cal O}_X$,
       are Morita equivalent.
    \end{itemize}
 \end{itemize}

\bigskip

\noindent $\bullet$
Re-examine Example~1.5.
 %
 \begin{itemize}
  \item[$\cdot$]
   Any existing way in noncommutative algebraic geometry
    to define the topological space $\Space(M_r({\Bbb C}))$
    for the ring $M_r({\Bbb C})$
   implies that $\Space(M_r({\Bbb C}))=\{\pt\}=\Spec{\Bbb C}$,
   if one really wants to define $\Space(M_r({\Bbb C}^r))$ honestly.

  \item[$\cdot$]
   One is thus supposed to define
    a morphism from the ringed space $(\Spec{\Bbb C}, M_r({\Bbb C}))$
    to $({\Bbb A}^1_{\Bbb C},{\cal O}_{{\Bbb A}^1_{\Bbb C}})$
    by a pair $(f,f^{\sharp})$,
   where
    $f:\Spec{\Bbb C}\rightarrow {\Bbb A}^1_{\Bbb C}=\Spec({\Bbb C}[y])$
     and
    $f^{\sharp}: {\cal O}_{{\Bbb A}^1_{\Bbb C}}
                  \rightarrow f_{\ast}(M_r(\Bbb C))$.

  \item[$\cdot$]
   Since $f_{\ast}(M_r({\Bbb C}))$ is a skyscraper sheaf at $f(\pt)$,
   the data $(f,f^{\sharp})$ is the same
    as the data of a ${\Bbb C}$-algebra homomorphism
    $$
     h:{\Bbb C}[y]\rightarrow M_r({\Bbb C})
    $$
    such that $\Ker h=h^{-1}(0)\subset {\Bbb C}[y]$
     is the ideal associated to a fuzzy point supported
     at $f(\pt)\in {\Bbb A}^1_{\Bbb C}$.
   This is a subclass of morphisms in Example~1.5
    which assume the additional constraint that $I_{\varphi}=((y-c)^n)$
    for some $c\in {\Bbb C}$ and $1\le n\le r$.

  \item[$\cdot$]
   Mathematically, there is nothing wrong with this.$\,/\!/\hspace{1ex}$
   But,
    for our purpose even just to describe D0-branes
     on the complex line ${\Bbb A}^1_{\Bbb C}$,
    this is {\it too restrictive}.$\,/\!/\hspace{1ex}$
   In particular,
    we won't be able to reproduce the Higgsing/un-Higgsing
     nor the D-brane recombination phenomenon
    if we confine ourselves to this traditional definition
     of morphisms between ringed spaces.
 \end{itemize}

\bigskip

\noindent $\bullet$
Morphisms between ringed spaces: second attempt guided by D-branes.
 \begin{itemize}
  \item[$\cdot$]
   {\it Forget}(!) the topological space;
   {\it keep only the rings}.

  \item[$\cdot$]
   A ``{\it morphism}"
    $\varphi:(X, {\cal O}_X)\rightarrow (Y,{\cal O}_Y)$
    is {\it defined contravariantly by}
    a ``{\it morphism}"
     $\varphi^{\sharp}:{\cal O}_Y\rightarrow {\cal O}_X$  {\it
    in the sense of an equivalence class of gluing systems
     of ring-homomorphisms}, {\sl when the latter can be defined}.

  \item[$\cdot$]
   In the commutative case, this recovers the usual definition
    of morphisms between (commutative) schemes since
   in that case $\varphi^{\sharp}$, in the sense above, truly defines
   a compatible continuous map (with respect to the Zariski topology)
   $\varphi:X\rightarrow Y$ {\it and} a sheaf homomorphism
     ${\cal O}_Y\rightarrow \varphi_{\ast}{\cal O}_X$,
     the usual $\varphi^{\sharp}$ in the theory of (commutative) schemes.
 \end{itemize}

\bigskip

\noindent $\bullet$ A major issue: {\it localization of an
(associative, unital) noncommutative ring.}
 \begin{itemize}
  \item[$\cdot$]
   We are thinking of a `space', whatever that means, {\it contravariantly}
   as an equivalence class of {\it gluing systems of rings} related by
   localizations of rings.

  \item[$\cdot$]
   An unfortunate fact:
   The notion of {\it localization of an (associative, unital)
   noncommutative ring} begins in 1931 in a work of Ore
   and is much more subtle than in the commutative case.

  \item[$\cdot$]
   Various techniques were developed, e.g.\ Gabriel's filter construction.
   This is an ongoing issue for the current and the future ring-theorists.
 \end{itemize}

\bigskip

\noindent $\bullet$
A D-brane-sensible/motivated settlement
 in the inperfect noncommutative world:\\
re-reading Proto-Definition~1.4.
 %
 \begin{itemize}
  \item[$\cdot$]
   Keep track only of and glue rings only through
    {\it central localizations};\\
   i.e.\ localizations only by elements that are in the center of a ring.

  \item[$\cdot$]
   $(X,{\cal O}_X^{nc})$,
    where $X$ is a topological space
          with a commutative structure sheaf ${\cal O}_X$\\
          $\mbox{$\hspace{4.2em}$}$
          that lies in the center of ${\cal O}_X^{nc}$,\\
    $\mbox{$\hspace{1em}$}=$
     an equivalence class of gluing system of rings
     in which the localization uses elements in ${\cal O}_X$.

  \item[$\cdot$]
   The topological space $X$ is only auxiliary and for this purpose.\\
   Truly, we are thinking the space $\Space({\cal O}_X^{nc})$,
    though we never define it!\\
   This explains basic noncommutative geometry on the D-brane world-volume.

  \item[$\cdot$]
   For the target-space-time $Y$, take any class of commutative or
    noncommutative spaces as long as they have a presentation
    as a class of gluing system of rings.

  \item[$\cdot$]
   A morphism $(X,{\cal O}_X^{nc})\rightarrow Y$
    is defined contravariantly as an equivalence class
     of gluing systems of ring-homomorphisms,
    exactly as one does for schemes.
 \end{itemize}

\bigskip

\noindent $\bullet$
A shift of perspective: a comparison with functor of points:
 \begin{itemize}
  \item[$\cdot$]
   In commutative algebraic geometry, we are very used to
    the concept that a space can also be defined
    by how others spaces are mapped into it.$\,/\!/\hspace{1ex}$
   Here, we are taking a reverse perspective.
   As indicated by Example~1.5, we are actually using how a ``space"
    can be mapped to other (more understood) spaces to feel
    this hidden-behind-the-veil ``space".
 \end{itemize}

\bigskip

\begin{flushleft}
{\bf Reflection and a conjcture on D-geometry in the sense of Douglas:\\
     Douglas meeting Polchinski-Grothendieck.}
\end{flushleft}
Before leaving this section,
 for the conceptual completeness of the lecture,
 let me give also some reflection on the notion of
 `{\it D-geometry}' in the sense of Michael Douglas [Do].
For any $r\in {\Bbb N}$,
 this is meant to be a certain noncommutative K\"{a}hler geometry
 on the moduli/configuration space ${\cal X}_r$ of D-brane
 for $r$-many D-branes on a K\"{a}hler manifold;
see [Do] and [D-K-O] for a more detailed description.
Let me recall first some basic facts from
 [L-L-L-Y] (D(2)) and [L-Y7] (D(6)).

\begin{lemma}
{\bf [special role of D0-brane moduli stack].}
{\rm ([L-L-L-Y: Sec.~3.1] (D(2)) and [L-Y7: Sec.~2.2] (D(6)).)}
 Let
  $Y$ be a (commutative) scheme over ${\Bbb C}$ and
  ${\frak M}^{{0^{A\!z}}^f}_r\!\!(Y)$
   be the moduli stack of D0-branes of rank/type $r$ on $Y$
   in the sense of Proto-Definition~1.4.
 Then,
  a morphism
   $$
    \varphi\; :\;
     (X,{\cal O}_X^{A\!z},{\Bbb C}^r)\; \longrightarrow\; Y\,,
   $$
   as defined in Proto-Definition~1.4
  is specified by a morphism
   $$
    \tilde{\varphi}\; :\;
     X\; \longrightarrow\; {\frak M}^{{0^{A\!z}}^f}_r\!\!(Y)\,;
   $$
   and vice versa.
\end{lemma}

Note that
 the universal family of D0-branes on $Y$
   over ${\frak M}^{{0^{A\!z}}^f}_r\!\!(Y)$
  defines an Azumaya structure sheaf
   ${\cal O}_{{\frak M}^{{0^{A\!z}}^f}_r\!\!(Y)}$
  with a fundamental module ${\cal E}_{{\frak M}^{{0^{A\!z}}^f}_r\!\!(Y)}$
  on ${\frak M}^{{0^{A\!z}}^f}_r\!\!(Y)$,
 realizing it canonically
  as an Azumaya (Artin/algebraic) stack with a fundamental module.
A comparison of
  the space-time aspect -- cf.~Aspect (2) in {\sc Figure}~1-3,
  the setting of [Do] and [D-K-O],  and
  the above lemma
 leads one then to the following conjecture,
  which brings Douglas' D-geometry into our setting:

\begin{conjecture}
{\bf [D-geometry: Douglas meeting Polchinski-Grothendieck].}
 An atlas for the Azumaya stack with a fundamental module
  $$
   (\,{\frak M}^{{0^{A\!z}}^f}_r\!\!(Y)\,,\,
      {\cal O}_{{\frak M}^{{0^{A\!z}}^f}_r\!\!(Y)}\,,\,
      {\cal E}_{{\frak M}^{{0^{A\!z}}^f}_r\!\!(Y)}\,)\; :=\; Y^{nc}_r
  $$
  corresponds to the configuration space ${\cal X}_r$ of D-branes
   in the work of Douglas [Do].
 For $Y$ a K\"{a}hler manifold,
 there exists an associated formal K\"{a}hler geometry
  on the irreducible component of $Y^{nc}_r$
  that contains all the $0$-dimensional ${\cal O}_Y$-module of length $r$
  whose support are $r$ distinct points on $Y$.
 This associated formal K\"{a}hler geometry can be made
    to satisfy the mass conditions of [Do] and [D-K-O]
  if and only if the K\"{a}hler manifold $Y$ is Ricci flat.
\end{conjecture}

\bigskip

\section{Azumaya geometry as the origin of the master nature\\ of D-branes.}

\noindent $\bullet$
In Sec.~1, we see that the matrix/Azumaya-type noncommutativity
 on D-brane world-volume occur in a very fundamental
  - almost the lowest - level.$\,/\!/\hspace{1ex}$
We also see in Example~1.5 that
 thinking of D-branes on an open-string target-space-time $Y$
 as morphisms from such Azumaya-type noncommutative space
 with a fundamental module does reproduce some features of D-branes
 in string theory.

\bigskip

\noindent $\bullet$
If the setting is truly correct from string-theory point of view,
 we should be able to see what string-theorists see
 in quantum-field-theory language solely by our formulation.
In particular,
 \begin{itemize}
  \item[$\cdot$] {\bf Q.}\
  \parbox[t]{13cm}{\it {\bf [QFT vs.\ maps]}\\
   Can we reconstruct the geometric object that arises
    in a quantum-field-theoretical study of D-branes
    through morphisms from Azumaya noncommutative spaces?}
 \end{itemize}
This is the guiding question for this section.

\bigskip

\begin{flushleft}
{\bf Azumaya noncommutative geometry
     as the origin of the master nature of D-branes.}
\end{flushleft}
$\bullet$
During the decade I was struggling to understand D-branes,
I read through quite a few string-theorists's work
 with various level of understanding.
However, there is one thing I failed to come by at that time:
 \begin{itemize}
  \item[{\bf Q.}] {\it
   For those D-brane works that carry a strong flavor of geometry,
   what exactly is going on\\ geometrically?}
 \end{itemize}
For that reason, for the scattered small pieces about D-branes
 I felt I understood something, I remained missing a real crucial piece
 to link them.
For that reason, I didn't truly understand what D-brane really is.
I asked several string theorists,
 including {\it Joe Polchinski} in TASI 1996,
 {\it Jeffrey Harvey} in TASI 1999,
 {\it Ashoke Sen} in TASI 2003,
 {\it Paul Aspinwall}'s TASI 2003 lectures and after-lecture
 discussions with participants, and {\it Cumrun Vafa}
 in a few occasions in and outside his courses at Harvard.
Each one gave me an answer.
That means each of these experts has his own working definition
 of D-branes strong and encompassing enough to create lots of
 significant works.
Yet, I wasn't able to fit their answer coherently together
 even to the picture I obtained
 when I read these experts' work.$\,/\!/\hspace{1ex}$
Then came a completely unexpected twist in the end of 2006.
A train of communications with {\it Duiliu-Emanuel Diaconescu}
 on a vanishing lemma of open Gromov-Witten invariants
  derived from [L-Y1] and [L-Y2]
 and his joint work with Florea [D-F] on open-string world-sheet instantons
 from the large $N$ duality of compact Calabi-Yau threefolds
 drove me back to re-understand D-branes.
After leaving this project for four years,
in this another attempt I came up with the understanding that
 there is a very fundamental noncommutativity on the D-brane world-volume
 and D-branes can be thought of as morphisms from such spaces,
 if this notion of morphism is defined ``correctly".
Then, I re-looked at some of the works
 that influenced me but I had failed to understand the true geometry behind.
At last, these pieces settle down coherently by one single notion:
 namely, {\it morphisms from Azumaya spaces}$\,$!

Below are a few examples.

\bigskip

\noindent $\bullet$
{\it For B-branes}$\,$: (Cf.\ [L-Y7: Sec.~2.4] (D(6)).)

\bigskip

\noindent
{\bf (1)}
\parbox[t]{15cm}{{\bf
 Bershadsky-Sadov-Vafa:
 Classical and quantum moduli space of D$0$-branes.}\\
 ({\it Bershadsky-Sadov-Vafa vs.\ Polchinski-Grothendieck}$\,$;
 [B-V-S1], [B-V-S2], [Va] (1995).)}
 \begin{itemize}
  \item[]
  The moduli stack ${\frak M}^{\,0^{Az^f}}_{\bullet}\!\!(Y)$
   of morphisms from Azumaya point with a fundamental module
   to a {\it smooth variety $Y$ of complex dimension $2$}
   contains various substacks with different coarse moduli space.
  One choice of such gives rise to the symmetric product
   $S^{\bullet}(Y)$ of $Y$
  while another choice gives rise to the Hilbert scheme
   $Y^{[\bullet]}$ of points on $Y$.
  The former play the role of the classical moduli and
   the latter quantum moduli space of D0-branes studied
  in [Va] and in [B-V-S1], [B-V-S2].

  \vspace{-1.6ex}
  \item[] $\hspace{1.2em}$
  See [L-Y3: Sec.~4.4] (D(1)),
       theme: `A comparison with the moduli problem of gas of D0-branes
       in [Va] of Vafa' for more discussions.
 \end{itemize}

\bigskip

\noindent
{\bf (2)}
\parbox[t]{15cm}{{\bf
 Douglas-Moore and Johnson-Myers:\\
 D-brane probe to an ADE surface singularity.}\\
({\it Douglas-Moore/Johnson-Myers
      vs.\ Polchisnki-Grothendiecek}$\,$;
 [D-M] (1996), [J-M] (1996).)}
 \begin{itemize}
  \item[]
   Here, we are compared with the setting of Douglas-Moore [Do-M].
   The notion of `morphisms from an Azumaya scheme
    with a fundamental module'
    can be formulated as well when the target $Y$ is a stack.
   In the current case, $Y$ is the {\it orbifold associated to an
    ADE surface singularity}. It is a {\it smooth Deligne-Mumford stack}.
   Again, the stack ${\frak M}^{\,0^{Az^f}}_{\bullet}\!\!(Y)$
    of morphisms from Azumaya points with a fundamental
    module to the orbifold $Y$ contains various substacks
    with different coarse moduli space.
   An appropriate choice of such gives rise to the resolution of
    ADE surface singularity.

  \vspace{-1.6ex}
  \item[] $\hspace{1.2em}$
  See [L-Y4] (D(3)) for a brief highlight of [D-M],
   details of the Azumaya geometry involved, and more references.
  In Sec.~3 of this lecture, we will present an abundance conjecture
   extracted and purified from the study initiated by [D-M].
 \end{itemize}

\bigskip

\noindent
{\bf (3)}
\parbox[t]{15cm}{{\bf
 Klebanov-Strassler-Witten: D-brane probe to a conifold.}\\
({\it Klebanov-Strassler-Witten vs.\ Polchinski-Grothendieck}$\,$;
 [K-W] (1998), [K-S] (2000).)}
 \begin{itemize}
  \item[]
   Here, the problem is related to the moduli stack
    ${\frak M}^{\,0^{Az^f}}_{\bullet}\!\!(Y)$
    of morphisms from Azumaya points with a fundamental module
    to a local conifold $Y$, a singularity Calabi-Yau $3$-fold,
    whose complex structure is given by
    $Y=\Spec({\Bbb C}[z_1,z_2,z_3,z_4]/(z_1z_2-z_3z_4))$.
   Again, different resolutions of the conifold singularity of $Y$
    can be obtained by choices of substacks from
    ${\frak M}^{\,0^{Az^f}}_{\bullet}\!\!(Y)$,
   as in Tests (1) and (2).
   Such a resolution corresponds to a low-energy effective geometry
    ``observed" by a stacked D-brane probe to $Y$
    when there are no fractional/trapped brane sitting
    at the singularity {\boldmath $0$} of $Y$.

  \vspace{-1.6ex}
  \item[] $\hspace{1.2em}$
   New phenomenon arises when there are fractional/trapped  D-branes
   sitting at {\boldmath $0$}. Instead of resolutions of the conifold
   singularity of $Y$, a low-energy effective geometry ``observed"
   by a D-brane probe is a complex deformation of $Y$ with topology
   $T^{\ast}S^3$ (the cotangent bundle of $3$-sphere).
   From the Azumaya geometry point of view, two things happen:
   \begin{itemize}
    \item[$\cdot$]
     Taking {\it both} the (stacked-or-not) D-brane probe
       and the trapped brane(s) into account,
     the Azumaya geometry on the D-brane world-volume remains.

    \item[$\cdot$]
     A noncommutative-geometric enhancement of $Y$ occurs
      via morphisms
      $$
      \xymatrix{
       & \Xi=\Space R_{\Xi} \ar[d]^{\pi^{\Xi}} \\
       Y \ar @{^{(}->}[r] & {\Bbb A}^4  &.
      }
      $$
   \end{itemize}
   Here,
    ${\Bbb A}^4=\Spec({\Bbb C}[z_1,z_2,z_3,z_4])$,
    $$
     R_{\Xi}\;=\;
      \frac{ {\Bbb C}\langle\,\xi_1,\xi_2,\xi_3,\xi_4\, \rangle }
      { ([\xi_1\xi_3, \xi_2\xi_4]\,,\, [\xi_1\xi_3, \xi_1\xi_4]\,,\,
         [\xi_1\xi_3, \xi_2\xi_3]\,,\, [\xi_2\xi_4, \xi_1\xi_4]\,,\,
         [\xi_2\xi_4, \xi_2\xi_3]\,,\, [\xi_1\xi_4, \xi_2\xi_3]) }
    $$
     with
      ${\Bbb C}\langle\,\xi_1,\xi_2,\xi_3,\xi_4\, \rangle$
       being the associative (unital) ${\Bbb C}$-algebra
       generated by $\xi_1,\,\xi_2,\,\xi_3,\,\xi_4$ and
      $[\bullet\,,\bullet^{\prime}\,]$ being the commutator,
     $Y\hookrightarrow {\Bbb A}^4$ via the definition of $Y$ above,
       and
     $\pi^{\Xi}$ is specified by the ${\Bbb C}$-algebra homomorphism
      $$
      \begin{array}{cccccl}
       \pi^{\Xi,\sharp}  & :
        & {\Bbb C}[z_1,z_2,z_3,z_4]  & \longrightarrow  &  R_{\Xi}\\
       && z_1 & \longmapsto  & \xi_1\xi_3 \\
       && z_2 & \longmapsto  & \xi_2\xi_4 \\
       && z_3 & \longmapsto  & \xi_1\xi_4 \\
       && z_4 & \longmapsto  & \xi_2\xi_3  & .\\
      \end{array}
      $$
   One is thus promoted to studying the stack
    ${\frak M}^{\,0^{Az^f}}_{\bullet}\!\!(\Space R_{\Xi})$,
    of morphisms from Azumaya points with a fundamental module to
    $\Space R_{\Xi}$.

  \vspace{-1.6ex}
  \item[] $\hspace{1.2em}$
   To proceed, we need the following notion:

  \item[]
  {\bf Definition 2.3.1.
  [superficially infinitesimal deformation].} {\rm
   Given associative (unital) rings,
    $R=\langle\,r_1,\,\ldots\,,r_m\,\rangle/\!\!\sim$ and $S$,
   that are finitely-presentable and
    a ring-homomorphism $h:R\rightarrow S$.
   A {\it superficially infinitesimal deformation} of $h$
     {\it with respect to the generators}
     $\{r_1,\,\ldots\,,r_m\}$ {\it of} $R$
    is a ring-homomorphism $h_{\varepsilon}:R\rightarrow S$
     such that
      $h_{\varepsilon}(r_i)=h(r_i)+\varepsilon_i$
       with $\varepsilon_i^2=0$,
      for $i=1,\,\ldots\,,m$.
    } 

  \item[]
   When $S$ is commutative,
    a superficially infinitesimal deformation of
     $\, h_{\varepsilon}:R\rightarrow S\,$
     is an infinitesimal deformation of $h$
     in the sense that $h_{\varepsilon}(r)=h(r)+\varepsilon_r$
      with $(\varepsilon_r)^2=0$, for all $r\in R$.
   This is no longer true for general noncommutative $S$.
   The $S$ plays the role of the Azumaya algebra
    $M_{\bullet}({\Bbb C})$ in our current test.
   It turns out that
   a morphism $\varphi:\pt^{A\!z}\rightarrow \Space R_{\Xi}$
    that projects by $\pi^{\Xi}$ to the conifold singularity
    {\boldmath $0$}$\in Y$ can have
     superficially infinitesimal deformations $\varphi^{\prime}$
    such that the image $(\pi^{\Xi}\circ \varphi^{\prime})(\pt^{A\!z})$
     contains not only {\boldmath $0$} but also points in
     ${\Bbb A}^4-Y$.
   Indeed there are abundant such superficially infinitesimal
    deformations.
   Thus, beginning with a substack ${\cal Y}$ of
    ${\frak M}^{\,0^{Az^f}}_{\bullet}\!\!(\Space R_{\Xi})$,
    that projects onto $Y$
    via $\varphi\mapsto \Image(\pi^{\Xi}\circ\varphi)$,
    one could use a $1$-parameter family
    of superficially infinitesimal deformations of $\varphi\in {\cal Y}$
    to drive ${\cal Y}$ to a new substack ${\cal Y}^{\prime}$
    that projects to {\boldmath $0$}$\cup Y^{\prime}\subset {\Bbb A}^4$,
    where $Y^{\prime}$ is smooth (i.e.\ a deformed conifold).
   It is in this way that a deformed conifold $Y^{\prime}$
    is detected by the D-brane probe via the Azumaya structure
    on the common world-volume of the probe and the trapped brane(s).

  \vspace{-1.6ex}
  \item[] $\hspace{1.2em}$
  See [L-Y5] (D(4)) for a brief highlight of [K-W] and [K-S],
   details of the Azumaya geometry involved, and more references.
 \end{itemize}

\bigskip

\noindent
{\bf (4)}
\parbox[t]{15cm}{{\bf
 G\'{o}mez-Sharpe:
 Information-preserving geometry, schemes, and D-branes.}\\
({\it G\'{o}mez-Sharpe vs.\ Polchisnki-Grothendieck}$\,$;
 [G-S] (2000).)}
 \begin{itemize}
  \item[]
   Among the various groups who studied the foundation of D-branes,
    this is a work that is very close to us in spirit.
   There, G\'{o}mez and Sharpe began with the quest:
    [G-S: Sec.~1]
   \begin{quote} ``{\it
    As is well-known, on $N$ coincident D-branes,
     $U(1)$ gauge symmetries are enhanced to $U(N)$ gauge symmetries,
    and scalars that formerly described normal motions of the branes
     become $U(N)$ adjoints.
    People have often asked what the deep reason for this behavior is
     -- what does this tell us about the geometry seen by D-branes?
     }",
   \end{quote}
   like us.
   They observed by comparing colliding D-branes with
    colliding torsion sheaves in algebraic geometry
    that it is very probable that
    \begin{quote}
     {\it coincident D-branes
     should carry some fuzzy structure} --
     {\it perhaps a nonreduced scheme structure}
    \end{quote}
    though the latter may carry more information
    than D-branes do physically.
   Further study on such nilpotent structure was done in [D-K-S];
   cf.~[L-Y7: Sec.~4.2: theme
    `The generically filtered structure on the Chan-Patan bundle
     over a special Lagrangian cycle on
     a Calabi-Yau torus'] (D(6)).

  \vspace{-1.6ex}
  \item[] $\hspace{1.2em}$
  From our perspective,
   \begin{quote}
    {\it the (commutative) scheme/nilpotent structure
    G\'{o}mez and Sharpe proposed/ observed on a stacked D-brane
    is the manifestation/residual of the Azumaya (noncommutative)
    structure on an Azumaya space with a fundamental module
   when the latter forces itself into a commutative space/scheme
    via a morphism}.
   \end{quote}
  This connects our work to [G-S].
 \end{itemize}

\bigskip

\noindent
{\bf (5)}
\parbox[t]{15cm}{{\bf
 Sharpe: $B$-field, gerbes, and D-brane bundles.}\\
({\it Sharpe vs.\ Polchinski-Grothendieck}$\,$;
 [Sh] (2001).)}
 \begin{itemize}
  \item[]
   Recall that a $B$-field on the target space(-time) $Y$
    specifies a gerbe ${\cal Y}_B$ over $Y$ associated to
    an $\alpha_B\in \check{C}^2_{\et}(Y,{\cal O}_Y^{\ast})$
    determined by the $B$-field.
  A morphism $\varphi:(X^{\!A\!z},{\cal E})\rightarrow (Y,\alpha_B)$
   from a general Azumaya scheme with a twisted fundamental module
   to $(Y,\alpha_B)$ can be lifted to a morphism
   $\breve{\varphi}:({\cal X}^{A\!z},{\cal F})\rightarrow {\cal Y}_B$
   from an Azumaya ${\cal O}_X^{\ast}$-gerbe with a fundamental module
   to the gerbe ${\cal Y}_B$.
  In this way, our setting is linked to Sharpe's picture of
   gerbes and D-brane bundles in a $B$-field background.

  \vspace{-1.6ex}
  \item[] $\hspace{1.2em}$
  See [L-Y6: Sec.~2.2] (D(5)) theme:
  `The description in term of morphisms from Azumaya gerbes
    with a fundamental module to a target gerbe'
  for details of the construction.
 \end{itemize}

\bigskip

\noindent
{\bf (6)}
\parbox[t]{15cm}{{\bf
 Dijkgraaf-Hollands-Su{\l}kowski-Vafa: Quantum spectral curves.}\\
({\it Dijkgraaf-Hollands-Su{\l}kowski-Vafa
      vs.\ Polchisnki-Grothendieck}$\,$;}\\[.8ex]
                                         {\rm $\mbox{\hspace{3.6ex}}$
 [D-H-S-V] (2007), [D-H-S] (2008).)}
 \begin{itemize}
  \item[]
   Here we focus on a particular theme in these works: the notion of
    {\it quantum spectral curves from the viewpoint of D-branes}.
   Let
    $C$ be a smooth curve,
    ${\cal L}$ an invertible sheaf on $C$,
    ${\cal E}$ a coherent locally-free ${\cal O}_C$-module,  and
    {\boldmath ${\cal L}$}$=\boldSpec(\Sym^{\bullet\,}({\cal L}^{\vee}))$
     be the total space of ${\cal L}$.
   Here, ${\cal L}^{\vee}$ is the dual ${\cal O}_C$-module of ${\cal L}$.
   Then one has the following canonical one-to-one correspondence:
   $$
    \left\{
     \begin{array}{l}
      \mbox{${\cal O}_C$-module homomorphisms}\\
      \phi:{\cal E}\rightarrow {\cal E}\otimes{\cal L}
     \end{array}
    \right\}\;
     \longleftrightarrow\;
    \left\{
     \begin{array}{l}
       \mbox{morphisms
        $\varphi:(C^{A\!z},{\cal E})\rightarrow$
                                    {\boldmath ${\cal L}$}}\\
       \mbox{as spaces over $C$}
     \end{array}
    \right\}
   $$
   induced by the canonical isomorphisms
   $$
   \Hom_{{\cal O}_C}({\cal E},{\cal E}\otimes{\cal L})\;
     \simeq\; \Gamma({\cal E}^{\vee}\otimes{\cal E}\otimes{\cal L})\;
     \simeq\; \Hom_{{\cal O}_C}
               ({\cal L}^{\vee},\,\Endsheaf_{{\cal O}_C}({\cal E}))\,.
   $$
   Let $\Sigma_{({\cal E},\phi)}\subset\;${\boldmath ${\cal L}$}
    be the (classical) spectral curve
    associated to the Higgs/spectral pair $({\cal E},\phi)$;
   cf.~e.g.\ [B-N-R], [Hi], and [Ox].
   Then, for $\varphi$ corresponding to $\phi$,
    $\Image\varphi\subset \Sigma_{({\cal E},\phi)}$.
   Furthermore, if $\Sigma_{({\cal E},\phi)}$ is smooth,
    then $\Image\varphi=\Sigma_{({\cal E},\phi)}$.
   This gives a {\it morphism-from-Azumaya-space interpretation of
    spectral curves}.

  \vspace{-1.6ex}
  \item[] $\hspace{1.2em}$
   To address the notion of `quantum spectral curve',
    let ${\cal L}$ be the sheaf $\Omega_C$ of differentials on $C$.
   Then the total space {\boldmath $\Omega$}$_C$ of $\Omega_C$ admits
    a canonical ${\Bbb A}^1$-family
    $Q_{{\Bbb A}^1}\mbox{\boldmath $\Omega$}_C$
    of deformation quantizations with the central fiber
    $Q_0\mbox{\boldmath $\Omega$}_C=\mbox{\boldmath $\Omega$}_C$.
   Let
    $({\cal E},\phi:{\cal E}\rightarrow {\cal E}\otimes\Omega_C)$
     be a spectral pair  and
    $\varphi: (C^{A\!z},{\cal E})\rightarrow \mbox{\boldmath $\Omega$}_C$
     be the corresponding morphism.
   Denote the fiber of $Q_{{\Bbb A}^1}\mbox{\boldmath $\Omega$}_C$
    over $\lambda\in {\Bbb A}^1$ by
    $Q_{\lambda}\mbox{\boldmath $\Omega$}_C$.
   Then, due to the fact that the Weyl algebras are simple algebras,
    the spectral curve $\Sigma_{({\cal E},\phi)}$
    in {\boldmath $\Omega$}$_C$ in general may not have a direct
    deformation quantization into $Q_{\lambda}\mbox{\boldmath $\Omega$}_C$
    by the ideal sheaf of $\Sigma_{({\cal E},\phi)}$
    in ${\cal O}_{\mbox{\scriptsize\boldmath $\Omega$}_C}$
   since this will only give
    ${\cal O}_{Q_{\lambda}\mbox{\scriptsize\boldmath $\Omega$}_C}$,
     which corresponds to the empty subspace of
     $Q_{\lambda}\mbox{\boldmath $\Omega$}_C$.
   However, one can still construct
    an ${\Bbb A}^1$-family
     $(Q_{{\Bbb A}^1}C^{A\!z}, Q_{{\Bbb A}^1}{\cal E})$
     of Azumaya quantum curves with a fundamental module
     out of $(C^{A\!z},{\cal E})$ and
    a morphism
     $\varphi_{{\Bbb A}^1}:
      (Q_{{\Bbb A}^1}C^{A\!z}, Q_{{\Bbb A}^1}{\cal E})
      \rightarrow Q_{{\Bbb A}^1}\mbox{\boldmath $\Omega$}_C$
     as spaces over ${\Bbb A}^1$,
    using the notion of `$\lambda$-connections'
     and `$\lambda$-connection deformations of $\phi$',
    such that
    \begin{itemize}
     \item[$\cdot$]
      $\varphi_0:= \varphi_{{\Bbb A}^1}|_{\lambda=0}$
       is the composition
       $\,(Q_0C^{A\!z},Q_0{\cal E})\,
           \longrightarrow\, (C^{A\!z},{\cal E})
            \stackrel{\varphi}{\longrightarrow}
             \mbox{\boldmath $\Omega$}_C\,$,
       where\\
         $(Q_0C^{A\!z},Q_0{\cal E})\rightarrow (C^{A\!z},{\cal E})$
         is a built-in dominant morphism from the construction;

     \item[$\cdot$]
      $\varphi_\lambda := \varphi_{{\Bbb A}^1}|_{\lambda}\,
        :\: (Q_{\lambda}C^{A\!z},Q_{\lambda}{\cal E})
            \longrightarrow Q_{\lambda}\mbox{\boldmath $\Omega$}_C$,
      for $\lambda\in {\Bbb A}^1-\{\mathbf 0\}\,$,
      is a morphism of Azumaya quantum curves with a fundamental module
      to the deformation-quantized noncommutative space
      $Q_{\lambda}\mbox{\boldmath $\Omega$}_C$.
    \end{itemize}
   In other words, we {\it
    replace the notion of `quantum spectral curves' by
    `quantum deformation $\varphi_{\lambda}$ of the morphism
    $\varphi$'}.
   In this way, both notions of classical and quantum spectral curves
    are covered in the notion of morphisms from Azumaya spaces.

  \vspace{-1.6ex}
  \item[] $\hspace{1.2em}$
   See [L-Y6: Sec.~5.2] (D(5))
    for more general discussions, details, and more references.
 \end{itemize}

\bigskip

\noindent $\bullet$
{\it For A-branes}$\,$:

\bigskip

\noindent
{\bf (7)}
\parbox[t]{15cm}{{\bf
 Denef: (Dis)assembling of A-branes under a split attractor flow.}\\
({\it Denef-Joyce meeting\ Polchisnki-Grothendieck}$\,$;
 [De] (2001), [Joy1] (1999), [Joy2] (2002--2003).)}
 \begin{itemize}
  \item[]
   (Dis)assembling of A-branes under a split attractor flow
    is realizable as Morse cobordisms of morphisms from Azumaya spaces
    with a fundamental module into the family of Calabi-Yau $3$-folds
    associated to the flow in the complex moduli space of the Calabi-Yau.
   Cf.~[L-Y8: Sec.~3.2] (D(7)).
 \end{itemize}

\bigskip

\noindent
{\bf (8)}
\parbox[t]{15cm}{{\bf
 Cecotti-Cordava-Vafa: Recombination of A-branes under RG-flow.}\\
({\it Cecotti-Cordova-Vafa
      meeting Polchisnki-Grothendieck}$\,$;
 [C-C-V] (2011).)}
 \begin{itemize}
  \item[]
   The renormalization group flow (RG-flow) in their setting
    specifies a flow on the moduli stack of morphisms
    from an Azumaya $3$-sphere with a fundamental module
    to the Calabi-Yau $3$-folds in question.
   The associated deformation family of morphisms corresponds
    the their brane recombinations.
   Cf.~[L-Y8: Sec.~2.3] (D(7)), [L-Y9] (D(8.1)), and work in progress.
 \end{itemize}

\bigskip

These and many more examples together motivate the next theme.

\bigskip

\begin{flushleft}
{\bf Azumaya noncommutative algebraic geometry
     as the master geometry for commutative algebraic geometry.}
\end{flushleft}
$\bullet$
A surprising picture emerges:
 \begin{itemize}
  \item[$\cdot$]
   {\bf [unity in geometry vs.\ unity in string theory]}\\[1.6ex]
   \framebox[18.6em][c]{\parbox{17.6em}{\it
    the master nature of morphisms from\\ Azumaya-type
    noncommutative spaces\\ with a fundamental module in geometry}}
    \hspace{1em} in parallel to \hspace{1em}
   \framebox[9em][c]{\parbox{8em}{\it
    the master nature of D-branes in \\superstring theory}}
 \end{itemize}
This strongly suggests that
 \begin{itemize}
  \item[$\cdot$]
   {\it Azumaya noncommutative algebraic geometry could play the role
        as the master geometry for commutative algebraic geometry.}
 \end{itemize}
Details remain to be understood.


\bigskip

\section{D-brane resolution of singularities - an abundance conjecture.}

\begin{flushleft}
{\bf Beginning with Douglas and Moore:
     D-brane resolution of singularities.}
\end{flushleft}
$\bullet$
For this third part of the lecture,
 let me begin with the work of Douglas and Moore [D-M].
 \begin{itemize}
  \item[$\cdot$]
   Let $\Gamma\simeq {\Bbb Z}_r \subset SU(2)$ acting on ${\Bbb C}^2$,
    with the standard Calabi-Yau $2$-fold structure,
    by automorphisms in the standard way.
   Consider the open and closed string target-space-time
    of the product form ${\Bbb R}^{5+1}\times [{\Bbb C}^2/\Gamma]$ and
    an effective-space-time-filling D-brane world-volume supported
    by the locus ${\Bbb R}^{5+1}\times\mbox{\boldmath $0$}$,
    where {\boldmath $0$} is the singular point of ${\Bbb C}^2/\Gamma$.

  \item[$\cdot$]
   The action of the supersymmetric QFT on the D-brane world-volume
    has various sectors arising from both open and closed strings.
   It involves, among other multiplets,
    vector multiplets and hypermultiplets.

  \item[$\cdot$]
   The potential energy function $V$ of hypermultiplets
    can be obtain by integrating out the Fayet-Iliopoulos D-term
    in the vector multiplets from the action.
   The result involves scalar fields $\vec{\phi}_{\bullet}$
    from NS-NS twisted sectors.

  \item[$\cdot$]
   From this,
     by taking $V^{-1}(0)/\mbox{\it global symmtry}$,
   one obtains the moduli space ${\cal M}_{\vec{\zeta}_{\bullet}}$
    of D-brane ground states.
   It depends on the vacuum expectation value $\vec{\zeta}_{\bullet}$
    of the scalar fields $\vec{\phi_{\bullet}}$.

  \item[$\cdot$]
   For appropriate choices of $\vec{\zeta}_{\bullet}$,
    ${\cal M}_{\vec{\zeta}_{\bullet}}$
    gives a resolution of the singularity of ${\Bbb C}^2/\Gamma$.
 \end{itemize}

\bigskip

\begin{flushleft}
{\bf The richness and complexity of Azumaya noncommutative space.}
\end{flushleft}
$\bullet$
There are lots of contents hidden in the Azumaya cloud
 ${\cal O}_X^{A\!z}$ of an Azumaya space
 $(X,{\cal O}_X^{A\!z},{\cal E})$;
cf.\ {\sc Figure}~3-1.$\,/\!/\hspace{1ex}$
This is already revealed by how an Azumaya point $\pt^{A\!z}$
 can be mapped to other spaces in the sense of Proto-Definition~1.4
 and is the origin of D-brane resolution of singularities,
 from our point of view;
cf.~{\sc Figure}~3-2.

\bigskip

\begin{flushleft}
{\bf An abundance conjecture.}
\end{flushleft}
\begin{definition} {\bf [punctual D0-brane].} {\rm
 (Cf.~[L-Y10: Definition~1.4] (D(9.1)).)
 Let $Y$ be a variety over ${\Bbb C}$.
 By a {\it punctual} $0$-dimensional ${\cal O}_Y$-module,
 we mean a $0$-dimensional ${\cal O}_Y$-module ${\cal F}$
  whose $\Supp({\cal F})$ is a single point
  (with structure sheaf an Artin local ring).
 A {\it punctual D0-brane on $Y$} of rank $r$
  is a morphism $\varphi:(\Spec{\Bbb C},\End(E),E)\rightarrow Y$,
  where $E\simeq {\Bbb C}^r$,
  such that $\varphi_{\ast}E$ is a ($0$-dimensional)
  punctual ${\cal O}_Y$-module.
 Let ${\frak M}_r^{0^{A\!z^f}_{\;p}}\!\!(Y)$
  be the {\it stack of punctual D0-branes of rank $r$ on a variety $Y$}.
 It is an Artin stack with atlas constructed from Quot-schemes.
 There is a morphism
  $\pi_Y:{\frak M}^{0^{A\!z^f}_{\;p}}\!\!(Y) \rightarrow Y$
  that takes $\varphi$ to $\Supp(\varphi_{\ast}E)$
  with the reduced scheme structure.
 $\pi_Y$ is essentially the Hilbert-Chow/Quot-Chow morphism.
}\end{definition}

\bigskip

\noindent $\bullet$
In term of this,
 note that:
 \begin{itemize}
  \item[$\cdot$]
   {\it Looking only at the internal part,
   then each element in ${\cal M}_{\vec{\zeta}_{\bullet}}$
   corresponds to a punctual D0-brane on $[{\Bbb C}^2/\Gamma]\,$.}
 \end{itemize}
It follows that
 the result of Douglas and Moore [D-M] of D-brane resolution
  of ADE surface singularities reviewed above can be rephrased as:
 (resuming the notation ${\Bbb A}^2$ for the affine variety
  behind ${\Bbb C}^2$.)

\begin{proposition}
{\bf [Douglas-Moore: D-brane resolution of ADE singularities].}
 There is an embedding
   $\widetilde{{\Bbb A}^2/\Gamma}
    \rightarrow {\frak M}^{0^{A\!z^f}_{\;p}}_1([{\Bbb A}^2/\Gamma])$
  that descends to a resolution
   $\widetilde{{\Bbb A}^2/\Gamma}\rightarrow {\Bbb A}^2/\Gamma$
   of singularities of ${\Bbb A}^2/\Gamma$.
\end{proposition}

\bigskip

\noindent $\bullet$
This,
  together with
   other existing examples of D-brane resolution of singularities
    -- including the case of conifolds --
   and the richness and complexity of the stack
    ${\frak M}^{0^{A\!z^f}_{\;p}}_r(Y)$,
 motivates the following abundance conjecture:

\begin{conjecture} {\bf [abundance].}
 Let $Y$ be a reduced quasi-projective variety over ${\Bbb C}$.
 Then,
  any birational model $Y^{\prime}\rightarrow Y$ of and over $Y$
  factors through an embedding of $Y^{\prime}$
  into the moduli stack ${\frak M}^{0^{A\!z^f}_{\;p}}_r(Y)$
  of punctual D0-branes of rank $r$ on $Y$, for $r$ sufficiently large.
\end{conjecture}

In particular,

\begin{conjecture} {\bf [D0-brane resolution of singularity].}
 Let $Y$ be a reduced quasi-projective variety over ${\Bbb C}$.
 Then,
 any resolution $\rho: Y^{\prime}\rightarrow Y$
  of the singularities of $Y$ factors through an embedding
  of $Y^{\prime}$ into ${\frak M}^{0^{A\!z^f}_{\;p}}_r(Y)$,
  for $r$ sufficiently large.
\end{conjecture}

\bigskip

\noindent $\bullet$
As a simple test, one has the following proposition:

\begin{proposition}
 {\bf [D0-brane resolution of curve singularity].}
 {\rm ([L-Y10 (L-(Baosen Wu)-Yau): D(9.1), Proposition~2.1].)}
 Conjecture~3.4 holds in the case of curves over ${\Bbb C}$.
 Namely,
 let
  $C$ be a (proper, Noetherian) reduced singular curve over ${\Bbb C}$
  and
  $$
   \rho\;:\; C^{\prime}\; \longrightarrow\; C
  $$
  be the resolution of singularities of $C$.
 Then, there exists an $r_0\in {\Bbb N}$
  depending only on
   the tuple $(n_{p^{\prime}})_{\rho(p^{\prime})\in C_{sing}}$
    and a (possibly empty) set
    $\{
     \bii(p)\,:\, p\in C_{sing}\,,\,
                     \mbox{C has multiple branches at $p$}\,\}$,
   both associated to the germ of $C_{sing}$ in $C$,
 such that,
 for any $r\ge r_0$,
  there exists an embedding
  $\tilde{\rho}: C^{\prime}\hookrightarrow
   {\frak M}^{0^{A\!z^f}_{\;p}}_r\!\!(C)$
  that makes the following diagram commute:
  $$
   \xymatrix{
     &&  {\frak M}^{0^{A\!z^f}_{\;p}}_r\!\!(C) \ar[d]^-{\pi_C} \\
    C^{\prime}\hspace{1ex}\ar @{^{(}->}[rru]^-{\tilde{\rho}} \ar[rr]^-{\rho}
     && \hspace{1ex}C\hspace{1ex}  &.
   }
  $$
 Here,
 \begin{itemize}
  \item[$\cdot$] {\rm
   $n_{p^{\prime}}\in {\Bbb N}$, for $\rho(p^{\prime})\in$
    the singular locus $C_{sing}\subset C$, is a {\it multiplicity}
   related to how the graph $\Gamma_{\rho}$ of $\rho$ intersects
   $C^{\prime}\times\{\rho(p^{\prime})\}$ (scheme-theoretically)
   in the product $C^{\prime}\times C\,$;}

  \item[$\cdot$] {\rm
   $\bii(p)\in {\Bbb N}$ is the {\it branch intersection index}
   of $p\in C_{sing}$; it is the least upper bound
   of the length of the $0$-dimensional schemes
   from the (scheme-theoretical) intersections of pairs of
   distinct branches of $C$ at $p\,$.}
 \end{itemize}
\end{proposition}

\bigskip

\noindent $\bullet$
Two remarks I should mention:

\begin{remark} {\it $[\,$another aspect$\,]$.} {\rm
 (Cf.~[L-Y10: Remark~0.1] (D(9.1)).)
 It should be noted that there is another direction
  of D-brane resolutions of singularities (e.g.\ [As1], [Br], [Ch]),
  from the point of view of (hard/massive/solitonic) D-branes
  (or more precisely B-branes)
  as objects in the bounded derived category of coherent sheaves.
 Conceptually that aspect and ours (for which D-branes are soft
  in terms of string tension) are in different regimes
  of a refined Wilson's theory-space of $d=2$ supersymmetric field
  theory-with-boundary on the open-string world-sheet.
 Being so, there should be an interpolation between these two aspects.
 It would be very interesting to understand such details.
}\end{remark}

\begin{remark}
{$[\,$string-theoretical remark$\,]$.} {\rm
 (Cf.~[L-Y10: Remark~1.7] (D(9.1)).)
 A standard setting (cf.\ [D-M]) in D-brane resolution of singularities
  of a (complex) variety $Y$ (which is a singular Calabi-Yau space
  in the context of string theory) is to consider
   a super-string target-space-time
    of the form ${\Bbb R}^{(9-2d)+1}\times Y$  and
   an (effective-space-time-filling) D$(9-2d)$-brane
    whose world-volume sits in the target space-time
    as a submanifold of the form ${\Bbb R}^{(9-2d)+1}\times\{p\}$.
 Here, $d$ is the complex dimension of the variety $Y$
  and $p\in Y$ is an isolated singularity of $Y$.
 When considering only the geometry of the internal part of this setting,
  one sees only a D0-brane on $Y$.
 This explains the role of D0-branes in the statement of
  Conjecture~1.5 and Conjecture~1.6.
 On the physics side, the exact dimension of the D-brane
  (rather than just the internal part) matters
 since
  supersymmetries and their superfield representations
   in different dimensions are not the same  and, hence,
  dimension does play a role in writing down
   a supersymmetric quantum-field-theory action for the world-volume
   of the D$(9-2d)$-brane probe.
 In the above mathematical abstraction, these data are now reflected
  into the richness, complexity, and a master nature of the stack
  ${\frak M}^{0^{A\!z^f}_{\;p}}_r\!\!(Y)$
   that is intrinsically associated to the internal geometry.
 The precise dimension of the D-brane as an object sitting in
  or mapped to the whole space-time becomes irrelevant.
}\end{remark}

\vspace{6em}

\begin{flushleft}
{\bf Epilogue.}
\end{flushleft}
In view of the fundamental role of Azumaya geometry for D-branes
 and the fact that Azumaya noncommutativity is lost under Morita equivalence
 and for that reason, most standard noncommutative algebraic geometers
 current days who follow the categorical language
 don't treat it as a significant noncommutative geometry,
one cannot help making the following moral,
 derived from Lao-Tzu (600 B.C.), {\sl Tao-te Ching}
   ({\sl The Scripture on the Way and its Virtue}), Chapter~11:

\bigskip

\centerline{\it What's naught could be the most useful!}

\bigskip

\vspace{16em}

\newpage

$ $

\vspace{16em}

\begin{figure}[htbp]
 \epsfig{figure=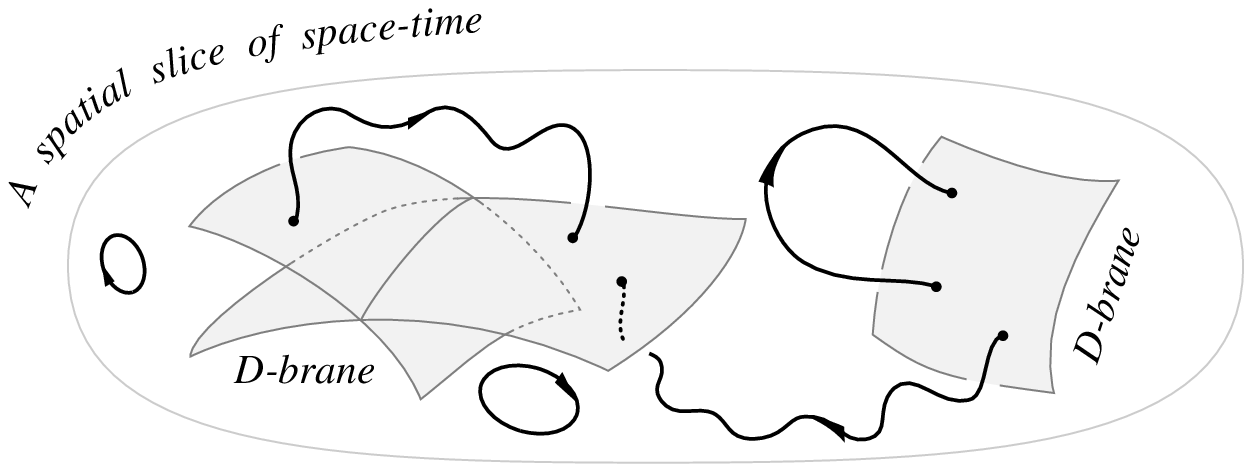,width=16cm}
 \centerline{\parbox{13cm}{\small\baselineskip 11pt
  {\sc Figure} 1-1.
  D-branes as boundary conditions for open strings in space-time.
  This gives rise to interactions of D-brane world-volumes
   with both open strings and closed strings.
  Properties of D-branes,
    including the quantum field theory on their world-volume and
              deformations of such,
   are governed by open and closed strings via this interaction.
  Both oriented open (resp.\ closed) strings and
   a D-brane configuration are shown.
  }}
\end{figure}

\newpage

$ $

\vspace{6em}

\begin{figure}[htbp]
 \epsfig{figure=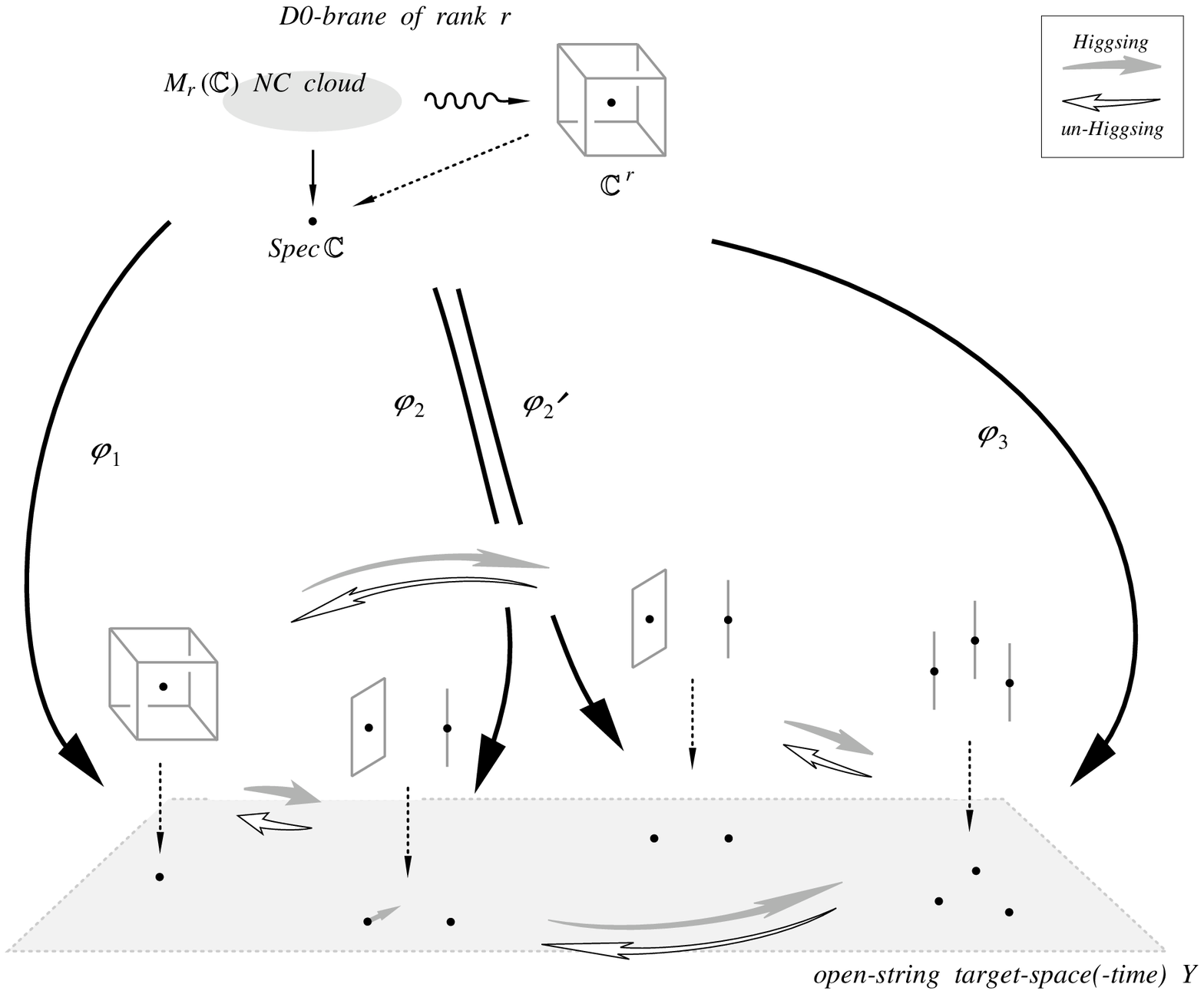,width=16cm}
 \centerline{\parbox{13cm}{\small\baselineskip 11pt
  {\sc Figure}~1-2. (Cf.\ [L-Y7: {\sc Figure}~2-1-1] (D(6)).)
 Despite that {\it Space}$\,M_r({\Bbb C})$ may look only
   one-point-like,
  under morphisms
  the Azumaya ``noncommutative cloud" $M_r({\Bbb C})$
  over {\it Space}$\,M_r({\Bbb C})$ can ``split and condense"
  to various image schemes with a rich geometry.
 The latter image schemes can even have more than one component.
 The Higgsing/un-Higgsing behavior of the Chan-Paton module of
   D$0$-branes on $Y$ ($={\Bbb A}^1$ in Example) occurs
  due to the fact that
   when a morphism
    $\varphi:$ {\it Space}$\,M_r({\Bbb C}) \rightarrow Y$
    deforms,
   the corresponding push-forward $\varphi_{\ast}E$
    of the fundamental module $E={\Bbb C}^r$
    on {\it Space}$\,M_r({\Bbb C})$ can also change/deform.
 These features generalize to morphisms
  from Azumaya schemes with a fundamental module to a scheme $Y$.
 Despite its simplicity, this example already hints at
  a richness of Azumaya-type noncommutative geometry.
 In the figure, a module over a scheme is indicated by
  a dotted arrow $\xymatrix{ \ar @{.>}[r] &}$.
 }}
\end{figure}

\begin{figure}[htbp]
 \epsfig{figure=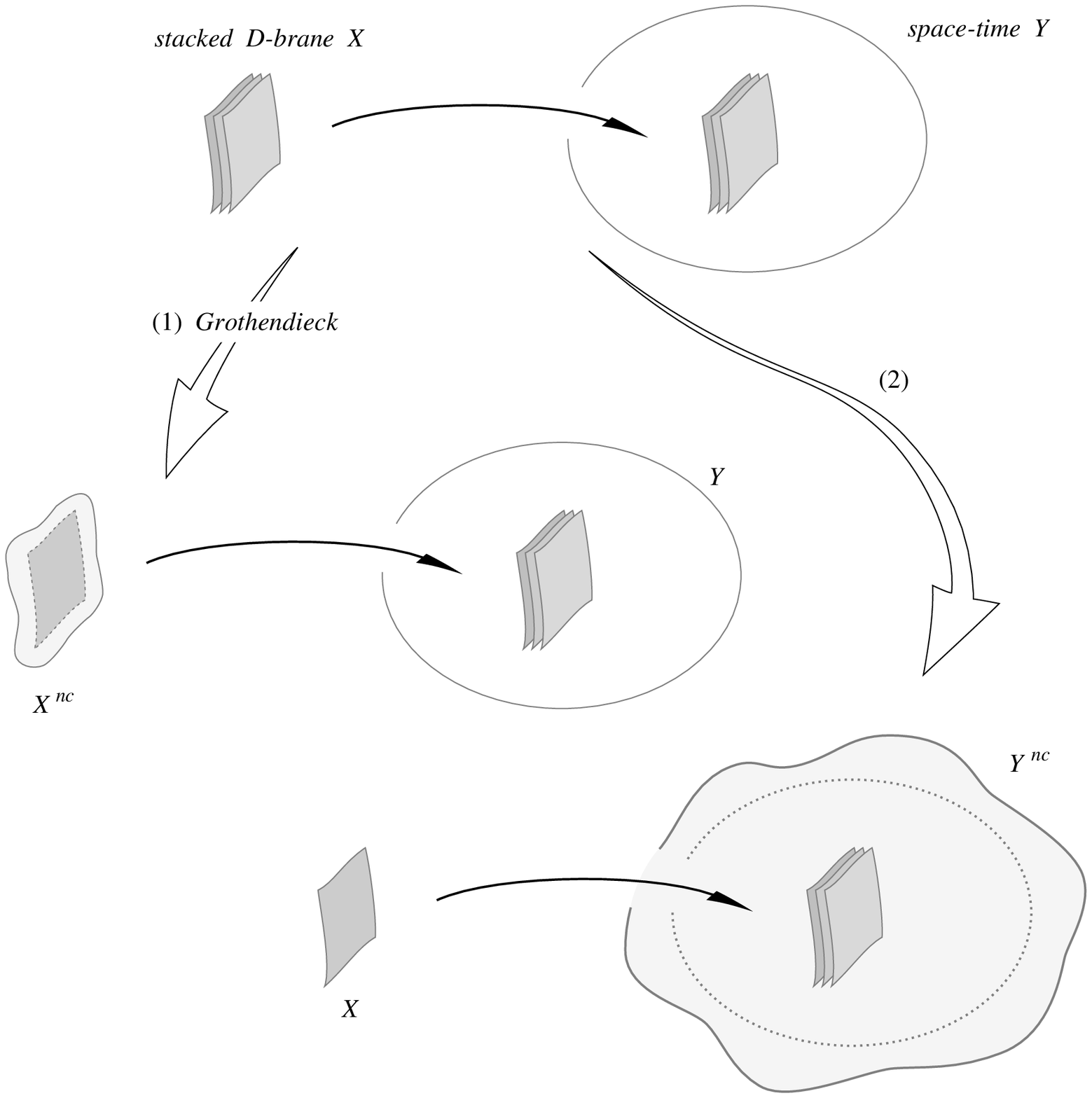,width=16cm}
 \centerline{\parbox{13cm}{\small\baselineskip 11pt
  {\sc Figure}~1-3. (Cf.\ [L-Y7: {\sc Figure}~1-1-2] (D(6)).)
  Two counter (seemingly dual but not quite) aspects
   on noncommutativity related to coincident/stacked D-branes:
   (1) noncommutativity of D-brane world-volume
       as its fundamental/intrinsic nature
   versus
   (2) noncommutativity of space-time as probed by stacked D-branes.
  (1) leads to the {\it Polchinski-Grothendieck Ansatz} and
   is more fundamental from Grothendieck's viewpoint
   of contravariant equivalence of the category of local geometries
   and the category of function rings.
  The matrix/Azumaya structure on coincident D-brane world-volume
   was also found in the work of Pei-Ming Ho and Yong-Shi Wu
   [P-W] (1996) in their own path.
  Their significant observation was unfortunately ignored
    by the majority of string-theory community.
  The latter pursued Path (2), following a few equally pival works
   including [Do] (1997) of Michael Douglas.}}
\end{figure}

\begin{figure}[htbp]
 \epsfig{figure=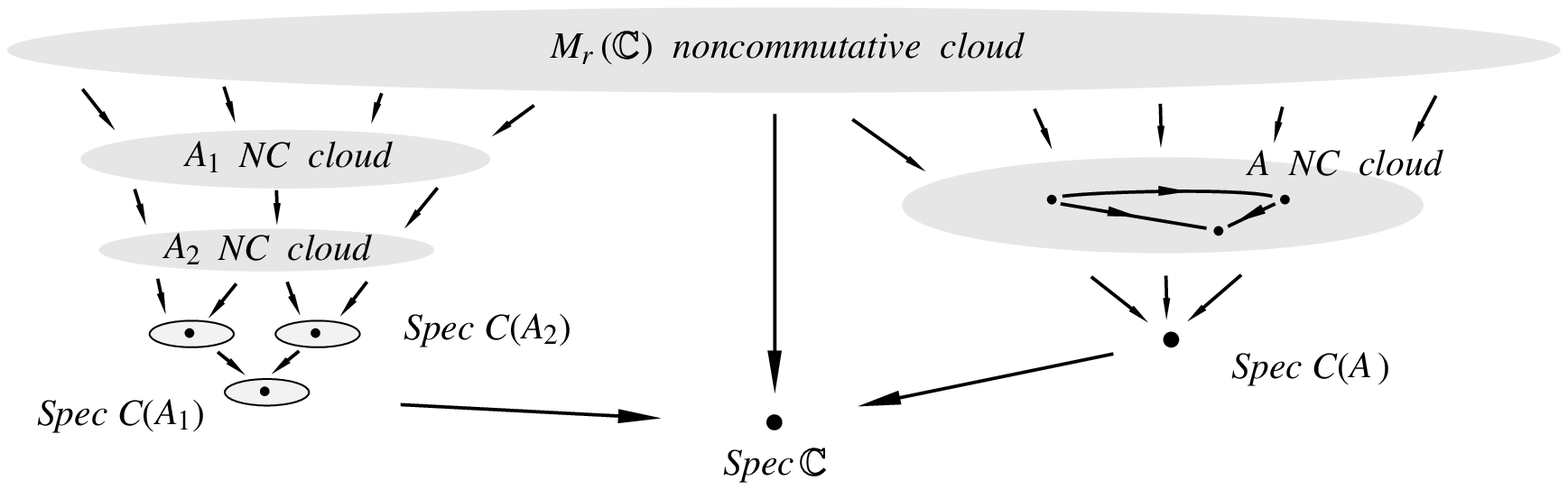,width=16cm}

 \vspace{3em}

 \centerline{\parbox{13cm}{\small\baselineskip 11pt
  {\sc Figure}~3-1. (Cf.~[L-Y4: {\sc Figure}~0-1] (D(3)).)
  An Azumaya scheme contains a very rich amount of geometry,
   revealed via its surrogates;
   cf.\ [L-L-S-Y: {\sc Figure}~1-3].
  Indicated here is the geometry of an Azumaya point
   $\pt^{\Azscriptsize} := (\smallSpec{\Bbb C}, M_r({\Bbb C}))$.
  Here, $A_i$ are ${\Bbb C}$-subalgebras of $M_r({\Bbb C})$
    and $C(A_i)$ is the center of $A_i$ with
   $$
     \begin{array}{cccccccl}
      M_r({\Bbb C}) & \supset  & A_1  & \supset  &  A_2
                    &\supset   & \cdots \\
       \cup  && \cup && \cup \\
     {\Bbb C}\cdot {\mathbf 1} & \subset  & C(A_1)
                    & \subset  & C(A_2)   & \subset  & \cdots & .
     \end{array}
   $$
  According to the Polchinski-Grothendieck Ansatz,
   a D$0$-brane can be modelled prototypically
   by an Azumaya point with a fundamental module of type $r$,
    $(\smallSpec{\Bbb C},\smallEnd({\Bbb C}^r),{\Bbb C}^r)$.
  When the target space $Y$ is commutative,
   the surrogates involved are commutative ${\Bbb C}$-sub-algebras
    of the matrix algebra $M_r({\Bbb C})=\End({\Bbb C}^r)$.
  This part already contains an equal amount of
   information/richness/complexity
   as the moduli space of $0$-dimensional coherent sheaves
   of length $r$.
  When the target space is noncommutative,
   more surrogates to the Azumaya point will be involved.
  Allowing $r$ to go to $\infty$ enables Azumaya points to probe
   ``infinitesimally nearby points" to points on a scheme
   to arbitrary level/order/depth.
  In (commutative) algebraic geometry,
   a resolution of a scheme $Y$ comes from a blow-up.
  In other words,
   a resolution of a singularity $p$ of $Y$ is achieved
   by adding an appropriate family of
   infinitesimally nearby points to $p$.
  Since D-branes with an Azumaya-type structure
   are able to ``see" these infinitesimally nearby points
    via morphisms therefrom to $Y$,
   they can be used to resolve singularities of $Y$.
  Thus, from the viewpoint of Polchinski-Grothendieck Ansatz,
   the Azumaya-type structure on D-branes is
   why D-branes have the power to ``see" a singularity
   of a scheme not just as a point,
   but rather as a partial or complete resolution of it.
  Such effect should be regarded as a generalization of
   the standard technique in algebraic geometry
   of probing a singularity of a scheme by arcs
   of the form $\Spec({\Bbb C}[\varepsilon]/(\varepsilon^r))$,
   which leads to the notion of jet-schemes
   in the study of singularity and birational geometry.}}
\end{figure}

\begin{figure}[htbp]
 \epsfig{figure=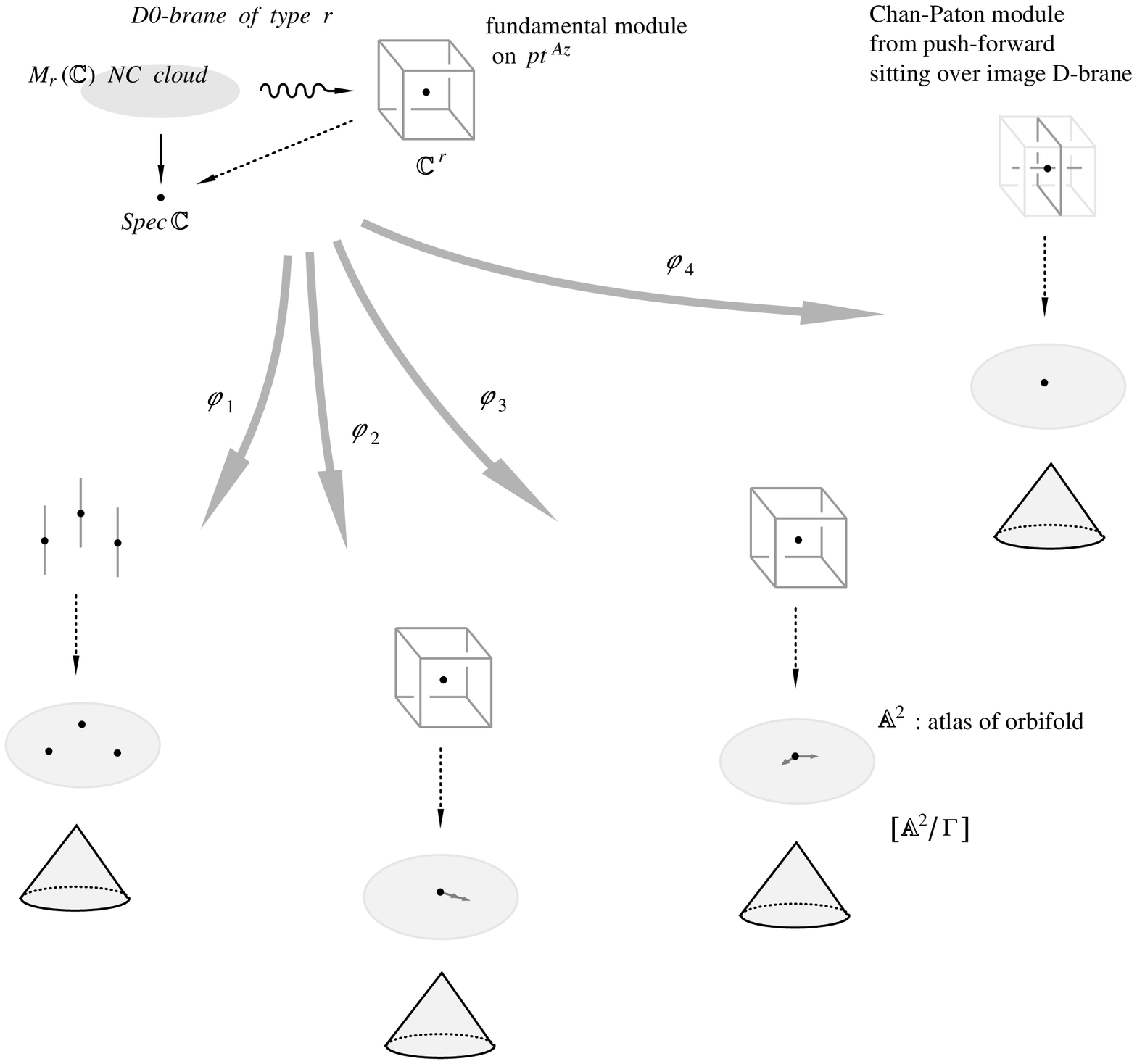,width=16cm}
 \centerline{\parbox{13cm}{\small\baselineskip 11pt
  {\sc Figure}~3-2. (Cf.~[L-Y4: {\sc Figure}~2-1] (D(3)).)
  Examples of morphisms from an Azumaya point with a fundamental module
   $(\smallSpec{\Bbb C}, \smallEnd({\Bbb C}^r),{\Bbb C}^r)$,
    which models an intrinsic D$0$-brane
    according to the Polchinski-Grothendieck Ansatz,
   to the orbifold $[{\Bbb A}^2/\Gamma]$ are shown.
  Morphism $\varphi_1$ is in Case (a)
   while morphisms $\varphi_2$, $\varphi_3$, $\varphi_4$
    are in Case (b).
  The image D$0$-brane under $\varphi_i$ on the orbifold $[{\Bbb A}^2/\Gamma]$
    is represented by a $0$-dimensional $\Gamma$-subscheme of length $\le r$
    on the atlas ${\Bbb A}^2$ of $[{\Bbb A}^2/\Gamma]$.
  }}
\end{figure}

\newpage

\baselineskip 10pt
\begin{flushleft}
{\small\bf Notes and acknowledgements added after the workshop.}
\end{flushleft}
{\small
 This note was prepared before the lecture with only mild revision
  and addition after coming back to Boston.
 For that reason, it is intentionally kept lecture-like
  so that the readers can get to the key points and the key words
  immediately without being distracted by formality.
 When writing this note three days before the workshop,
  I had in mind of it as part of notes for a minicourse.
 For this particular workshop, I selected the main part of Sec.~1
  and quick highlight in Sec.~3 and presented them mainly on the blackboard
  so that the audience can think over and digest the concept in real time.
 A vote was cast after presenting very slowly Example~1.5 and Remark~1.6
  to decide whether the audience, particularly string-theorists,
  agree that my notion of D-branes following the line of Grothendieck does
   correctly reflect string-theorists' D-branes (in the appropriate region
   of the related Wilson's theory-space, cf.\ beginning of Sec.~1).
 It turned out that there is no objection to the setting;
 yet it received only cautious acceptance:
 ``... can accept it but have to think more".
 This is another time I put the notion under the scrutinization
  of experts outside Yau's group and Harvard string-theory community
  since the first paper D(1) in the series that appeared in 2007.
 No objections do not necessarily imply believing it;
 there are still numerous themes in the series yet to be understood
  and completed.

 Special thanks to Charlie Beil for inviting me to this workshop,
  through which I learned many things I had been unaware of before;
 thanks also to many speakers who answer my various questions
  during or after their inspiring and resourceful lecture.
 Outside the workshop, I thank
  Paul Aspinwall
   for an illumination of a conceptual point in [As2]
    concerning central charge of B-branes;
  Ming-Tao Chuan
   for discussions on some technical issues on deformations
   of singular special Lagrangian cycles in Calabi-Yau manifolds
   related to D(8.1);
  Michael Douglas
   for illuminations/highlights of his D-geometry in [Do] and [D-K-O],
       explanation of a key question in [D-K-O]
         that requires a better understanding,
       and some reference guide
   -- indeed, though I am confident, it has been my wish to
   meet him directly to see if he has objections from physics ground
   to what I have been doing on D-branes;
  Pei-Ming Ho and Richard Szabo
   for preview of the note before the workshop;
  David Morrison
   for a discussion on some conceptual points
    on supersymmetric quantum field theory and Wilson's theory-space;
  Eric Sharpe
   for communicating to me a train of insights/comments on resolutions
   of singularities in string theory related to D(9.1)
   after I emailed him a preliminary version of this note
   before the workshop;  and
  Paul Smith for correcting my ridiculously wrong picture
   of the history of noncommutative algebraic geometry
   through and after his lecture and a literature guide
   -- there are clearly many things I have yet to learn.

 Comments/corrections/objections to this preliminary lecture note may
  be sent to the following as part of the basis for its future
  revision/improvement (after the project is pushed far enough):
 \\ \\

 \noindent
 {\sc e-mail}$\,:$
  chienliu@math.harvard.edu,$\;$
  chienhao.liu@gmail.com
} 

\newpage
\baselineskip 12pt

\end{document}